\documentclass[11pt]{article}
\usepackage[a4paper, margin=1in]{geometry}
\usepackage{bm,graphicx,amsmath,amsfonts,hyperref,amssymb,mathrsfs,tikz-cd,dsfont,amsthm,float,tikzit}

\tikzstyle{box}=[fill=white, draw=black, shape=rectangle]
\tikzstyle{longbox}=[fill=white, draw=black, shape=rectangle, minimum width=6cm]
\tikzstyle{3cmbox}=[fill=white, draw=black, shape=rectangle, minimum width=3cm]
\tikzstyle{2cmbox}=[fill=white, draw=black, shape=rectangle, minimum width=2cm]
\tikzstyle{4cmbox}=[fill=white, draw=black, shape=rectangle, minimum width=4cm]

\tikzstyle{dashy}=[-, densely dashed, tikzit draw=blue]
\tikzstyle{white}=[-, line width=4pt, draw=white, tikzit draw={rgb,255: red,191; green,255; blue,0}]
\tikzstyle{dotty}=[-, densely dotted, tikzit draw=red]
\tikzstyle{red}=[-, draw=red]
\tikzstyle{arrow}=[<-]

\usepackage{tikz}
\usetikzlibrary{knots}
\usepackage[eprint=true, url=false,doi=false,isbn=false]{biblatex}
\DeclareMathOperator{\Rep}{Rep}
\DeclareMathOperator{\PaB}{PaB}
\DeclareMathOperator{\wPaB}{wPaB}
\DeclareMathOperator{\CoB}{CoB}
\DeclareMathOperator{\wCoB}{wCoB}

\newcommand{\Z}{\mathbb{Z}}

\newcommand{\Hom}{\text{Hom}}

\newcommand{\B}{\mathcal{B}}

\numberwithin{equation}{section}
\newtheorem{theorem}{Theorem}[section]
\newtheorem{definition}{Definition}[section]

\newtheorem{remark}[theorem]{Remark}
\newtheorem{proposition}[theorem]{Proposition}
\newtheorem{lemma}[theorem]{Lemma}
\newtheorem{corollary}[theorem]{Corollary}
\newtheorem*{proof*}{Proof}

\newtheorem{manualtheoreminner}{Theorem}
\newenvironment{manualtheorem}[1]{%
	\renewcommand\themanualtheoreminner{#1}%
	\manualtheoreminner
}{\endmanualtheoreminner}

\title{Weak braiding for algebras in braided monoidal categories}
\author{Devon Stockall}

\begin{document}
\begin{center}
	\vspace*{1cm}
	
	\textbf{Weak braiding for algebras in braided monoidal categories}
	
	\vspace{0.5cm}
	
	Devon Stockall\footnote{stockall@imada.sdu.dk}
	
\end{center}
	
	\abstract{
	Under appropriate conditions, if one picks a commutative algebra $A$ with action of group $G$ in braided monoidal category $\mathcal{C}$, the category of $A$ modules in $\mathcal{C}$ obtains a natural crossed G-braided structure.  In the case of general commutative algebra object $A$ in braided monoidal category $\mathcal{C}$, one might ask what weakened notion of braiding one obtains on the category of $A$ modules in $\mathcal{C}$, and the relation that this braiding has to categorical symmetries acting on the associated quantum field theories, and on the algebra $A$ itself. 
	
	 In the following article, we present a definition of weak-braided monoidal category.  It is proven that braided $G$-crossed categories, and categories of modules over commutative algebra objects in braided monoidal categories are weak-braided monoidal categories.  Further, it is proven that, under reasonable assumptions, if a weak-braided category $\mathcal{D}$ is given by twisting the braiding by a collection of `twisting functors', then $\mathcal{D}$ is crossed-braided by a group.}
	\newpage
	\section{Introduction}
	Braided monoidal categories arise naturally in the description of symmetries of 2 and 3 dimensional quantum field theories.  
	
	It has been known since \cite{MooreSeiberg} that vertex operators of a 2-dimensional conformal field theory should arrange themselves in a braided monoidal category.  From a formal mathematical perspective, one associates to a 2d CFT a vertex operator algebra (VOA) V, and finds that, under appropriate finiteness restrictions, a category $\mathcal{C}$ of $V$ modules is a braided monoidal category  \cite{HL1}.
	
	An analogous story holds for 3-dimensional topological quantum field theory, developing first in Chern-Simons theory \cite{QFTandJonesPoly} and its formalization \cite{ReshetikhinTuraev}.  One finds that line operators in a 3d TQFT arrange themselves into a braided monoidal category, with objects corresponding to line operators of the theory, tensor product given by fusion of lines, and morphisms given by point operators at junctions of line operators (See also \cite{Fredenhagen1989,Gabbiani1990} and \cite{Kitaev2006} appendix E). 
	
	In both the two and three dimensional cases, it is natural to consider commutative algebra object $A$ in the corresponding braided monoidal category $\mathcal{C}$.  One then constructs the category $\Rep_\mathcal{C}A$ of $A$ modules in $\mathcal{C}$.  In the 2d CFT case, the algebra $A$ corresponds to an extension of the symmetry algebra of the conformal symmetry, and $\Rep_{\mathcal{C}}A$ corresponds again to labels of (possibly twisted) field insertions in the CFT with extended symmetry \cite{HKL,SVOAExtension}.  In the 3d TQFT case, the algebra $A$ corresponds to a family of anyons, and $\Rep_\mathcal{C}A$ corresponds to operators in the 3d TQFT obtained by condensing these anyons \cite{Kong2013}.   
	
	One defines the subcategory $\Rep^{loc}_\mathcal{C}A$, consisting of modules which braid trivially with $A$, ie $(M,\mu_M)\in \Rep_\mathcal{C}A$ satisfying
	\begin{equation}
		\mu_M\mathcal{M}_{A,M}=\mu_M
	\end{equation}
	Where $\mathcal{M}_{A,M}=\mathcal{R}_{M,A}\circ \mathcal{R}_{A,M}$ is the categorical monodromy.  In the 3d TQFT case, $\Rep^{loc}_\mathcal{C}A$ consists of line operators in the condensed theory, while in the 2d CFT interpretation $Rep^{loc}_\mathcal{C}A$ is the subcategory of $\Rep_\mathcal{C}A$ of local modules.  Only the subcategory $\Rep^{loc}_\mathcal{C}$ inherits the structure of a braided monoidal category from $\mathcal{C}$.  We argue that all objects in the category $\Rep_\mathcal{C}A$ are physical should inherit some weakened notion of braiding. 
	
	 To do this, we turn our attention to rational 2d conformal field theories, and their associated 3d TQFTs.  Viewing the Verlinde lines of \cite{Verlinde1988} as line operators in a bulk 3d TQFT, and the Kramers-Wannier duality defect of \cite{Schweigert2004} as a line operator in the 3d TQFT interpolating between the chiral and antichiral halves of the Ising model, one sees that line operators labeled by a module over the symmetry algebra should also be considered as (potentially \textit{non-invertible}) symmetry operators.
	
	 We next turn to the special case that the algebra object $A\in \mathcal{C}$ has the action of finite group $G$, whose fixed points are the monoidal unit $\mathds{1}\in \mathcal{C}$, as is the case in $G$-equivariant TQFTs and orbifold CFTs \cite{Dijkgraaf1989,Dijkgraaf1991,MTCandOrbifold, MTCandOrbifold2}.  Under appropriate assumptions, the simple objects $M\in \Rep_\mathcal{C} A$ are \underline{G-twisted} in the sense that
	\begin{equation}\label{Gtwist}
		\mu_M(g\boxtimes 1_M)\mathcal{M}_{A,M}=\mu_M
	\end{equation}
	for some $g\in G$.  This motivates the name `twisted module' for all non-local objects of $\Rep_\mathcal{C}A$ in the general case.  One also finds that $\Rep_\mathcal{C}A$ is a \underline{braided G-crossed} category \cite{McRae_2020} in the sense of \cite{HFTandGCrossed}, so that $\Rep_\mathcal{C}A$ inherits a `braiding up to G action' from $\mathcal{C}$. 
	
	In \cite{GeneralizedGlobalSymmetries}, it is suggested that all topological defects in a QFT should be considered as a sort of generalized symmetry.  In particular, usual grouplike symmetries should be associated to defects of codimension 1.  Contrasting the grouplike and non-invertible symmetry examples, we see that Kramers-Wannier and Verlinde defect lines, labeled by modules, are associated to codimension 2 defects in a 3d TQFT, while the grouplike defects of orbifolds are associated to codimension 1 defects in the associated 3d TQFT.  In particular, a $g$-twisted module, for $g\in G$, can be viewed as a line operator bounding an invertible surface operator which acts by the group element $g$.   
	
	Returning to the general case, we argue that all objects $\Rep_\mathcal{C}A$ should be seen as lines: codimension 2 objects in 3 dimensions, bounding some codimension 1 defect, and $\Rep_\mathcal{C}^{loc}A$ consists of lines bounding the trivial surface defect.  While the remaining twisted operators in $\Rep_\mathcal{C}A$ \textit{shouldn't} be braided in the usual sense, it is still physically reasonable to move one such codimension 2 defect operator around another, obtaining a braiding `up to' the action of the associated codimension 1 defect (which, may also act non-invertibly, leaving a tail).  Then $\Rep_\mathcal{C}A$ should also inherit some appropriately weakened notion of braiding structure.
	
	One would also like to study these generalized symmetries through the category $\Rep_\mathcal{C}A$.  One might also hope to find that the algebra object $A\in \mathcal{C}$ is acted on by some appropriate collection $K$ of non-invertible symmetries, and that 
	\begin{enumerate}
		\item Simple objects $M\in \mathcal{C}$ are $K$-twisted
		\item $\Rep_\mathcal{C}A$ is braided $K$-crossed
	\end{enumerate}
	It was proven in \cite{Andrew} that one naturally obtains a \underline{hypergroup} $K=\Rep_\mathcal{C}A//\Rep_\mathcal{C}^{loc}A$, and an action of $K$ on $A$.  While this does seem to be the appropriate collection of symmetries, one finds that the twisted modules $M\in \Rep_\mathcal{C}A$ are \textbf{not} twisted by this action in the sense of equation \ref{Gtwist}.  The action presented in \cite{Andrew} is the `invertible part' of a complete action of $K$ in $A$.
	
	The goal of the present work will be to address the means by which $\Rep_\mathcal{C}A$ is braided K-crossed, in the hopes that the sense in which modules are twisted by these symmetries appears more naturally there.  This should present itself as an action of the n-th braid group on n-th tensor products.  
	
	\subsection{Overview of the present work}
In section \ref{heuristic}, we argue why such a construction is necessary, and motivate the given structure.  In particular, we justify the use of operadic language for making this generalization of braiding.  In section \ref{operad}, we review the construction of the operads $\Omega_{\tilde{+}}$ and $\PaB_{\tilde{+}}$ of parenthesized permutations, and parenthesized braids, whose algebras are monoidal and braided monoidal categories, respectively.  We then construct an operad $\wPaB_{\tilde{+}}$, called the \underline{weak-braided monoidal category operad}, and then define a \underline{weak-braided} monoidal category to be an algebra over $\wPaB_{\tilde{+}}$.  For category $\mathcal{C}$, giving $\wPaB_{\tilde{+}}$ algebra stucture on $\mathcal{C}$ can be thought of as giving a coherent action of the $n$-th braid groups on `$n$-th braid-twisted tensor products'.
	
	In section \ref{weak braiding for extensions}, we develop categorical lemmata to simplify proofs in the vein of \cite{HKL,SVOAExtension}.  This will allow many of the results to follow from existence and universal properties of certain coequilizers.  These results are then used to prove the following theorem.
	
	\begin{manualtheorem}{\ref{Extensionweakbraiding}}
		Suppose that $\mathcal{C}$ is an abelian braided monoidal category with right exact tensor functors, and $A$ is an associative commutative algebra object in $\mathcal{C}$. Then $\Rep_\mathcal{C}A$ is weak-braided. 
	\end{manualtheorem}
	In section \ref{GCrossed}, the weak-braiding construction is used to address the special case of braided $G$-crossed fusion categories, as defined in \cite{HFTandGCrossed}.
	
	\begin{manualtheorem}{\ref{Gcrossedcase}}
		Suppose that $G$ is a group and that $\mathcal{D}$ is a braided $G$-crossed fusion category with additive structure functors.  Then $\mathcal{D}$ is weak-braided. 
	\end{manualtheorem}

The weak-braiding presented in this work may seem \textbf{overly} weak, and one might naively hope that, given a collection of non-invertible symmetries $K$, a braided $K$-crossed category would be given by an action of $K$ on category $\mathcal{C}$ by endofunctors, and $\mathcal{C}$ is `braided up to the action of $K$' as in the group case.  However, we prove the following partial converse to the previous statement.  
	\begin{manualtheorem}{\ref{Gcrossedsufficient}}
			Suppose that
		\begin{enumerate}
			\item $\mathcal{D}$ is linear, additive, semisimple and pivotal
			\item $\mathcal{D}$ is weak-braided with additive braid-tensor functors
			\item $S$ is a set of monoidal functors $F:\mathcal{D}\to \mathcal{D}$
			\begin{enumerate}
				\item For all $F\in S$ there exists simple $X\in \mathcal{D}$ such that 
				\begin{equation}
					(\nu,\mathcal{R})(\cdot,X)=(\nu,id)(F(\cdot),X)
				\end{equation}
				\item For all simple $X\in \mathcal{D}$, there is unique $F\in S$ such that 
				\begin{equation}
					(\nu,\mathcal{R})(\cdot,X)=(\nu,id)(F(\cdot),X)
				\end{equation}
				\item The tensor unit $\mathds{1}\in \mathcal{D}$ satisfies 
				\begin{equation}
					(\nu,\mathcal{R})(\cdot,\mathds{1})=(\nu,id)(id_{\mathcal{D}}(\cdot),\mathds{1})
				\end{equation}
			\end{enumerate}
		\end{enumerate}
		Then $S$ is a group and $\mathcal{D}$ is braided $S$-crossed.
	\end{manualtheorem}
	That is, if a weak-braiding structure on semisimple pivotal category $\mathcal{D}$ is exhibited by `twisting functors' $F:\mathcal{D}\to \mathcal{D}$, then there is group $G$ for which $\mathcal{D}$ is braided $G$-crossed.  Then if one hopes to study non-invertible symmetries $K$, then the category $\mathcal{C}$ can not be naively braided $K$-crossed.

	\subsection{Future questions}
	The present work leaves open the following questions:
	
	\begin{itemize}
		\item Given a weak-braided monoidal category $\mathcal{D}$, the associated non-invertible symmetries should be visible in `degeneracies' of the weak-braiding structure.  One might define a local subcategory $\mathcal{D}^{loc}$ of $\mathcal{D}$ to be the `kernel' of the weak braiding: a weak braided subcategory which is equivalent to a braided category.  Does $K=\mathcal{D}//\mathcal{D}^{loc}$ form a hypergroup, as in \cite{Andrew}?  
		\item Does the weak braiding structure on $\mathcal{D}$ illuminate how non-invertible symmetries \textbf{should} act on vertex algebras? 
		\item In the braided $G$-crossed case, one finds an `equivariantizion' procedure to recover $\mathcal{C}$ from $\Rep_{\mathcal{C}}A$.  Does such a procedure also exist for weak braiding?  
		\item One associates the operad of parenthesized braids $\PaB$ with configuration space of discs.  How does the weak braiding operad $\wPaB$ relate to a `configuration space' of universal covers of punctured discs?
		\item How is the present construction related to the pseudo-braided tensor categories of \cite{Soibelman1997}, especially in light of the work of \cite{Moriwaki2024}, which gives pseudo-braiding structure for categories of VOA modules which are $C_1$ cofinite with finitely generated duals?
		\item Non-invertible symmetries for algebra objects in 2 and 3 dimensions can likely be more naturally approached from a 3-categorical perspective, similar to \cite{CoreyPenneysReutter} and \cite{Kong2024}.  The collection of non-invertible symmetries labels objects in a monoidal two-category $\mathcal{K}$, and for symmetry $k\in \mathcal{K}$, k-twisted modules are morphisms $\mathds{1}\to k$, for $\mathds{1}\in \mathcal{K}$ the monoidal unit.
	\end{itemize}
	
	\subsection{Acknowledgements}
	I would like to thank Guillaume Laplante-Anfossi for helpful discussion.  Without him, the operad section of this paper would not be nearly so concise.  This work was supported by VILLUM FONDEN, VILLUM Young Investigator grant 42125. 
	
\section{A heuristic argument for weak braiding}\label{heuristic}
Here, we will demonstrate the necessity of the weak braiding construction in a special case.  The results proven later in this work will apply in much more generality.  We make the following technical assumptions only so that we have access to the graphical calculus of braid diagrams given in \cite{Turaev2016}, which lends itself well to heuristic argument. Suppose that $\mathcal{C}$ is a linear semisimple pivotal braided monoidal category and $A\in \mathcal{C}$ is a commutative algebra object.  

Recall that in braided monoidal category, one has the following form of `completeness relation':
	
	\ctikzfig{completeness}
	
	Where we sum over simple summands $N$ of $X\otimes Y$, and the morphisms $f_N:X\otimes Y\to M$ and $\tilde{f}_N:N\to X\otimes Y$ are left implicit for visual clarity.  Suppose that $X\in \Rep_{\mathcal{C}}A$ is simple.  Summing over simple summands $M$ of $X\otimes X^*$, we obtain a relation of the following form: 
	
	\scalebox{0.9}{\parbox{\linewidth}{%
	\ctikzfig{braidingcompleteness}}}
	
	Where we note that $A$ appears in this sum with multiplicity one, and use that $X$ is an $A$ module.  Notice the resemblance of this relation to the ones appearing \cite{MTCandOrbifold}.  In particular, when $X$ is $g$-twisted for g an automorphism of $A$, we see that the sum over $M\neq A$ vanishes.  
	
	Conversely, suppose that $X\in \Rep_{\mathcal{C}}A$ is such that this sum over $M\neq A$ vanishes which we will call the \underline{vanishing lasso condition}.  When the vanishing lasso condition holds, we also have the following:
	
	\ctikzfig{hopflinkhom}
	
	We see that, for module $X$ satisfying the vanishing lasso condition, the normalized `A open Hopf link' with $X$ acts as an $A$-homomorphism.  Assuming it does not vanish, it must act as an isomorphism, and in particular, the collection of A-open Hopf links given by such modules can be extended to a group in $Hom_\mathcal{C}(A,A)$.  
	
	We conclude that any module satisfying the vanishing lasso condition is `group twisted'.  In particular, any module `twisted by non-invertible symmetry' must remember the information of non-vanishing lassos, giving an explicit reason why twisted modules are not twisted by the action presented in \cite{Andrew} in the naive sense.  
	
	We now turn our attention to the structure of the associated categories: What features must a generalized `crossed-braiding' theory have to deal with non-vanishing lassos? 
	
	Consider $X,Y,Z\in \Rep_\mathcal{C}A$.  Whatever `braiding' we obtain on $\Rep_{\mathcal{C}}A$ should be induced from the braiding on $\mathcal{C}$, so the $\mathcal{C}$-morphism given by the following braid diagram should induce a $\Rep_{\mathcal{C}}A$ morphism from $X\otimes_A Y\otimes_A Z$ to an `appropriately twisted version' of the tensor product of $X,Y,Z$ over $A$.  We again use the completeness relation and argue that lassos should not vanish as above.
	
	\scalebox{0.9}{\parbox{\linewidth}{%
	\ctikzfig{Bparameternecessity}}}
	
	Then  we see that the necessary to `remember' information about $Z$ in order to recover the module map on $X$ in the twisted version of $X\otimes_AY\otimes_AZ$.  Since this can not be recovered as a morphism with additional lassos connecting only adjacent strands, one can not hope to recover this twisted version of $X\otimes_AY\otimes_AZ$ by iteratively tensoring together adjacent strands in a twisted way.  
	
	In this sense, these `braid twisted tensor products' of modules twisted by non-invertible symmetries are \underline{non-local}: one must remember information about non-adjacent strands.  This is essentially the statement that the braid operad is \textbf{not} generated as an operad by the crossing of two strands, even though it is generated as a group.

	One could then convince themselves that if they were to keep track of all over and under crossings of lassos, then the braid by which the strands are twisted could be recovered.  In fact, the braid itself encodes the information of crossings of lassos much more efficiently than one could hope to encode it directly.  This has the added benefit of applying also to cases that are not rigid or semisimple.  
	
	Then, with this and the non-locality of `twisted tensor products' in mind, when constructing a twisted tensor product for modules twisted by non-invertible symmetries, one should take as an input a family of $A$ modules $X_1,...,X_n\in \Rep_\mathcal{C}A$ and an element of the n-th braid group $\B\in B_n$.  One must then `twist' the multiplication on $X_1\otimes...\otimes X_n$ by the action of $\B$, and tensor all modules together simultaneously.  This construction should also have a coherent way to glue together and compose the associated braids. 

Abstracting away from the $\Rep_\mathcal{C}A$ example, one finds that the information of a weak-braided category should be as follows: 
\begin{itemize}
	\item A category $\mathcal{D}$
	\item Families of `twisted tensor functors' $\mathcal{D}^{\times n}\to \mathcal{D}$ parameterized by the $n$-th braid group $B_n$
	\item Rules for composition of twisted tensor functors, consistent with the parameterizing braids
\end{itemize}
This information is most naturally expressed by giving $\mathcal{D}$ the structure of an algebra over an appropriate operad.

\section{Operads and weak braiding}\label{operad}

It will be illustrative for us to construct the operad of (weakly unital) parenthesized braids $\PaB_{\tilde{+}}$ here, whose algebras are (weakly unital) braided monoidal categories, as to demonstrate exactly how the weak braiding operad $\wPaB_{\tilde{+}}$ differs.  

This is constructed approximately as follows:  We begin with the operad $\CoB_+$ in the category of groupoids, where the groupoids $\CoB_+(n)$ have permutations as objects, and braids interpolating between them as morphisms.  $\CoB_+$ algebras as strictly associative and unital braided monoidal categories.  We then consider an operad $\Omega_{\tilde{+}}$ in set, where the objects of $\Omega_{\tilde{+}}(n)$ are `parenthesized and weakly unital permutations'.  Pulling back the groupoid structure along the natural map $\tau:\Omega_{\tilde{+}}\to \CoB_+$ gives the operad $\PaB_{\tilde{+}}$ in the category of groupoids, whose algebras are (not necessarily strictly associative or strictly unital) braided monoidal categories.  

The construction of the weakly braided operad parallels this.  We begin with the operad $\wCoB_+$ in the category of groupoids, where the groupoids have \underline{colored braids} as objects, and braids interpolating between them as morphisms.  We then construct an operad $\Omega_{Br}$ in set, where the objects of $\Omega_{Br}(n)$ are parenthesized and weakly unital colored braids.  Pulling back the groupoid structure along the natural map $\pi:\Omega_{Br}\to \wCoB_+$ gives an operad $\wPaB_{\tilde{+}}$ in the category of groupoids, whose algebras are  (not necessarily strictly associative or strictly unital) \underline{weakly braided} monoidal categories. 

We will refer to \cite{Fresse2017} for operadic definitions.  A specific location for each definition is given here, for the benefit of those less familiar with operadic language.  
\begin{definition}\hfill
	\begin{itemize}
		\item We define a symmetric operad as in \cite{Fresse2017} 1.1.1.  The partial composition relations for operads, as are primarily used here, appear in 2.1.8.
		\item We define the symmetric group operad $\Sigma$ as given in \cite{Fresse2017} proposition 1.1.9.  
		\item For symmetric operad $P$, the symmetric group equivariance property takes the following form in terms of partial compositions in $\Sigma$ and $P$
		\begin{equation}
			\sigma_np_n\circ_i\sigma_mp_m=(\sigma_n\circ_i\sigma_m)(p_n\circ_{\sigma_n^{-1}(i)}p_m)
		\end{equation}
		for $p_n\in P(n)$, $p_m\in P(m)$, $\sigma_n\in \Sigma_n$, $\sigma_m\in \Sigma_m$, and $1\leq i\leq n$.
		\item For an object $\mathcal{C}$ in a symmetric monoidal category, there is a symmetric operad defined by $End_\mathcal{C}(n)=Hom(\mathcal{C}^{\otimes n}, \mathcal{C})$ with partial composition of morphisms.  An algebra over the operad $P$ is a morphism of operads $P\to End_\mathcal{C}$.
		\item We define braid groups $B_n$ on $n$ strands as in \cite{Fresse2017} 5.0.2, see also 5.0.9.  The pure braid group $P_n$ on $n$ strands is the kernel of this quotient map, and fits in short exact sequence of groups
		\begin{equation}\begin{tikzcd}
				1\arrow[r]&P_n\arrow[r,"\iota"]&B_n\arrow[r,"p_*"]&\Sigma_n\arrow[r]&1
		\end{tikzcd}\end{equation}
	The map $p_*$ assigns to each braid its \underline{underlying permutation}.  For $A\in B_n$, we will use the notation 
	\begin{equation}
		\overline{A}=p_*(\mathcal{B})
	\end{equation}
	\item We define the (non-symmetric) \underline{braid operad} $B_+$ as in \cite{Fresse2017} proposition 5.1.3.  The maps $p_*$ define a morphism of operads $B_+\to \Sigma_+$.  The partial composition in this operad satisfies the following equivariance property
	\begin{equation}
		A_nC_n\circ_iA_mC_m=(A_n\circ_iA_m)(C_n\circ_{\overline{A_n}^{-1}(i)}C_m)
	\end{equation}
\end{itemize}
\end{definition}

\subsection{Braided monoidal categories}

\begin{definition}
	For $n\in \Z_{\geq0}$, we consider the left action of $B_n$ on $\Sigma_n$ given by 
	\begin{equation}
		g\cdot u_n=u_np_*(g^{-1})
	\end{equation}
	For $u_n\in \Sigma_n$ and $g\in B_n$.  Then we define the groupoid of \underline{colored braids} on n strands $\CoB_+(n)$ to be the action groupoid of this action $\CoB_+(n)=\Sigma_n//B_n$.  That is, with objects and morphisms given by
	\begin{itemize} 
		\item $Ob(\CoB_+(n))=\Sigma_n$
		\item For $u,v\in \Sigma_n$ we have $Hom_{\CoB_+(n)}(u,v)=p_*^{-1}(v^{-1}u)$, with composition of morphisms given by multiplication in $B_n$.
	\end{itemize}
	Define the (unital) \underline{operad of coloured braids} $\CoB_+$ in the category of groupoidss by 
	\begin{itemize}
		\item Objects $\CoB_+(n)$ for $n\in Z_{\geq0 }$ 
		\item Symmetric group action, for $\sigma\in \Sigma_n$, given by left translation on objects, and on morphisms, is associated to the identity map 
		\begin{equation}
			Hom_{\CoB_+(n)}(u,v)=p_*^{-1}(v^{-1}u)=p_*^{-1}((\sigma v)^{-1}\sigma u)=Hom_{\CoB_+(n)}(\sigma u,\sigma v)
		\end{equation}
		\item Composition is given on objects by the composition in $\Sigma$, and on morphisms by\footnote{Notice that this differs slightly from \cite{Fresse2017}, page 184 as we believe there is a typo.  The composition on morphisms here is modified by the output permutation rather than the input.  This is easily checked to be consistent with other constructions in \cite{Fresse2017}.}
		\begin{multline}
			(g_n,g_m)\in Hom_{\CoB_+(n)}(u_n,v_n)\times Hom_{\CoB_+(m)}(u_m,v_m)\\
			\mapsto g_n\circ_{v_n^{-1}(i)}g_m\in Hom_{\CoB_+(n+m-1)}(u_n\circ_i u_m,v_n\circ_i v_m)
		\end{multline}
		Where $g_n\circ_{v_n^{-1}(i)} g_m$ is the partial composition in the braid group operad $B$.  
	\end{itemize}
	Note that $\overline{g_n\circ_ig_m}=\overline{g_n}\circ_i\overline{g_m}$ since the projection $p_*:B\to \Sigma$ is a morphism of operads.  This agrees with the definition given in \cite{Fresse2017} 5.2.8.
\end{definition}

\begin{definition}
	
	We define the \underline{weakly unital magma operad} $\Omega_{\tilde{+}}$ to be the symmetric operad freely generated by $\nu(x_1,x_2)\in \Omega_{\tilde{+}}(2),\ \nu(x_2,x_1)=(12)\cdot\nu(x_1,x_2)\in \Omega_{\tilde{+}}(2)$ and $x_0\in \Omega_{\tilde{+}}(0)$. 
\end{definition}
Note that this is a slight variation on \cite{Fresse2017} 6.1, which will allow us to have non-strict units in the resulting monoidal categories.  On the following, we use the notation $\nu=\nu(x_1,x_2)$.  

\begin{proposition}[\cite{Fresse2017} 1.2.4]
	If the base category $M$ has small colimits, then the category of symmetric operads in $M$ has small colimits. 
\end{proposition}

\begin{proposition}
	The category of groupoids is complete.
\end{proposition}

\begin{definition}
	Suppose $G,H$ are symmetric operads in the category of groupoids, and $\tau:G\to H$ is a morphism of symmetric operads.  Denote by $\tau^{-1}H$ the limit of the following diagram:
	\begin{equation}\begin{tikzcd}
			G\arrow[r,"\tau"]&H
	\end{tikzcd}\end{equation}
	In \cite{Fresse2017}, 6.1.5, $G$ is considered to be an operad in discrete groupoids, and this is called the pullback operad.  We notice that this is a categorical pullback square if it is instead written as
	\[\begin{tikzcd}
		\tau^{-1}H\arrow[r]\arrow[d]&H\arrow[d,"id"]\\
		G\arrow[r,"\tau"]&H
	\end{tikzcd}\]
\end{definition}

\begin{definition}
	 We may consider $\Omega_{\tilde{+}}$ as an operad in the category of groupoids by considering each $\Omega_{\tilde{+}}$ to be a discrete groupoid.  Since $\Omega_{\tilde{+}}$ is free, to define a morphism, it is sufficient to define it on generators.  
	
	We define morphism $\tau:\Omega_{\tilde{+}}\to \CoB_+$ of operads by 
	\begin{align}
		\tau_\mu=id_2\in \Sigma_2&& \tau_{x_0}=*\in \Sigma(0)
	\end{align}
	Where we note that evaluation is denoted by subscript.  We define the operad of \underline{weakly unital} \underline{parenthesized braids} $\PaB_{\tilde{+}}$ in the category of groupoids by 
	\begin{equation}
		\PaB_{\tilde{+}}=\tau^{-1}\CoB
	\end{equation}
	We note that $\Omega_{\tilde{+}}(n)$ and $\PaB_{\tilde{+}}(n)$ share the same objects, and for $p,q\in\PaB_{\tilde{+}}(n)$, the morphism sets are given by 
	\begin{equation}
		\Hom_{\PaB_{\tilde{+}(n)}}(p,q)=Hom_{\CoB_+(n)}(\tau(p),\tau(q))
	\end{equation}
\end{definition}

\begin{theorem}\cite{Yau_2021} Theorem 21.1.7.  
	Suppose that $\mathcal{C}$ is a small category.  Then a braided monoidal structure on $\mathcal{C}$ is equivalent to the structure of a $\PaB_{\tilde{+}}$ algebra on $\mathcal{C}$. 
\end{theorem}

\subsection{Weak braided monoidal categories}

\begin{definition}\label{wCoB}
	For $n\in \Z_{\geq 0}$, consider the set $\Sigma_n\times B_n$, with left action of $B_n$ given by
	\begin{equation}
		g\cdot (u,A)=(up_*(g^{-1}),gA)
	\end{equation}
	for $(u,A)\in \Sigma_n\times B_n$ and $g\in B_n$.  Define the groupoid $\wCoB_+(n)$ to be the action groupoid $\Sigma_n\times B_n//B_n$ by this action of $B_n$.  That is, defined by
\begin{itemize}
	\item $Ob(\wCoB_{+}(n))=\Sigma_n\times B_n$
	\item Morphisms given by 
	\begin{equation}
		Hom_{\wCoB_{+}}((u,A),(v,C))=p_*^{-1}(v^{-1}u)\cap \{CA^{-1}\}		
	\end{equation}
	for $(u,A),(v,C)\in \wCoB_{+}(n)$, with composition of morphisms given by multiplication in $B_n$. 
\end{itemize}
Define the operad $\wCoB_+$ by 
\begin{itemize}
	\item Symmetric group action defined on objects by $\sigma(u,A)=(\sigma u,A)$ for $(u,A)\in \wCoB_{\tilde{+}}(n)$ and $\sigma\in \Sigma_n$, and by the identity map 
	\begin{equation}
		p_*^{-1}(v^{-1}u)\cap\{CA^{-1}\}\to p_*^{-1}((\sigma v)^{-1}\sigma u)\cap \{CA^{-1}\}
	\end{equation}
	on morphisms. 
	\item Operadic compositions given on objects by
\begin{equation}
	(u_n,A_n)\circ_i(u_m,A_m)=(u_n\circ_iu_m,A_n\circ_{u_n^{-1}(i)}A_m)
\end{equation}
For $(u_n,A_n)\in \wCoB_+(n)$ and $(u_m,A_m)\in \wCoB^+(m)$. The operadic composition on morphisms is given by 
\begin{equation}
	(g_n,g_m)\mapsto g_n\circ_{v_n^{-1}(i)}g_m
\end{equation}   
which is the $v_n^{-1}(i)$-th partial composition in the braid operad. 
\item The unit functor is defined by 
\begin{equation}
	1\in \mathds{1}\mapsto (id_1,id_1)\in \Sigma(1)\times B_1 
\end{equation}
viewed as an object in $\wCoB_+(1)$.
\end{itemize}

\end{definition}

\begin{proof*}
	$\wCoB_+(n)$ is defined as an action groupoid, and so we need not check that it is a groupoid.  However, we note that, by the definition of action groupoid, we have 
	\begin{equation}
		g\in Hom_{\wCoB_+(n)}((u,A),(v,C))\implies gA=C
	\end{equation}
	resulting in the presentation given above.  
	
	We must then confirm that the composition structure given above defines a symmetric operad.  Notice $w\CoB_+$ is essentially an enhancement of $\CoB_+$ by braids, and so, many of the axioms here will follow from the axioms on $\CoB$ and $B$.

	We first note, using equivariance in the braid group operad, that 
	\begin{equation}
		id_{n+m-1}=A_nA_n^{-1}\circ_iA_mA_m^{-1}=(A_n\circ_iA_m)(A_n^{-1}\circ_{\overline{A_n}^{-1}(i)}A_m^{-1})
	\end{equation}
	so we conclude that $A_n^{-1}\circ_{\overline{A_n}(i)}A_m^{-1}=(A_n\circ_i A_m)^{-1}$.  Suppose $g\in Hom_{\wCoB_+(n)}((u_n,A_n),(v_n,C_n))$ and $h\in Hom_{\wCoB_+(m)}((u_m,A_m),(v_m,C_m))$ for $(u_n,A_n), (v_n,C_n)\in \wCoB_+(n)$ \\ and $(u_m,A_m),(v_m,C_m)\in \wCoB_+(m)$.  
	
	Note that $\overline{C_n}\ \overline{A_n}^{-1}=\overline{C_nA_n^{-1}}=v_n^{-1}u_n$ so that $u_n^{-1}=\overline{A_n}\ \overline{C_n}^{-1}v_n^{-1}$.  
	 
Then we have the following equality
	\begin{multline}
	g_n\circ_{v_n^{-1}(i)}g_m=C_nA_n^{-1}\circ_{v_n^{-1}(i)}C_mA_m^{-1}\\
	=(C_n\circ_{v_n^{-1}(i)}C_m)(A_n^{-1}\circ_{\overline{C_n}^{-1}v_n^{-1}(i)}A_m^{-1})
	=(C_n\circ_{v_n^{-1}(i)}C_m)(A_n\circ_{\overline{A_n}^{-1}\overline{C_n}v_n^{-1}(i)}A_m)^{-1}\\
	=(C_n\circ_{v_n^{-1}(i)}C_m)(A_n\circ_{u_n^{-1}(i)}A_m)^{-1}
	\end{multline}
	Then, since the usual coloured braid groupoids $\CoB_+$ define an operad, we have $g_n\circ_{v_n^{-1}(i)}g_m\in p_*^{-1}((v_n\circ_iv_m)^{-1}(u_n\circ_iu_m))$, and we conclude that 
	\begin{equation}
		g_n\circ_{v_n^{-1}(i)}g_m\in Hom_{\wCoB_+(n+m-1)}((u_n,A_n)\circ_i(u_m,A_m),(v_n,C_n)\circ_i(v_m,C_m))
	\end{equation}
	so the given composition is well defined.  Since the partial composition on morphisms is identical to in the usual colored braid operad $\CoB_+$, we conclude that this partial composition is functorial.  
	
	The proof of vertical and horizontal associativity, unit, and symmetric group equivariance is structurally identical, on the level of objects, to the proof given under definition \ref{omegabr}, and so we defer this portion of the proof to after definition \ref{omegabr}.  On the level of morphisms, we note that morphisms between any two objects in $\wCoB_{+}(n)$ are unique if they exist, and we have already proven that composition functors are well defined, so it follows also that vertical and horizontal associativity, unit, and symmetric equivariance axioms hold on the level of morphisms.  We conclude that $\wCoB_+$ is an operad in the category of groupoids.  
\end{proof*}

\begin{definition}
	 $p\in \Omega_{\tilde{+}}(n)$ may be written uniquely as 
	\begin{equation}
		p=p^{len}(x_{t_p(1)},...,x_{t_p(len(p))})
	\end{equation}
	where 
	\begin{itemize}
		\item $p^{len}\in \Omega_{\tilde{+}}(len(p))$ is written only using compositions of $\nu$. 
		\item $t_p:\{1,2,...,len(p)\}\to \{0,1,...,n\}$ is bijective except on $t_p^{-1}\{0\}$.  
	\end{itemize}
	Abusing notation, we define 
	\begin{equation}
		t_p^{-1}(i)=t_p|_{\{1,...,len(p)\}-t_p^{-1}\{0\}}^{-1}(i)
	\end{equation}
	for $i=1,...,n$.  Then $len(p)$ counts the number of characters occurring in the word $p$, and $t_p^{-1}(i)$ gives the location of the character $x_i$ in the word $p$.  
\end{definition}

\begin{definition}\label{omegabr}
	For $n\in \Z_{\geq0}$, define the operad $\Omega_{Br}$ in $Set$ by
	\begin{itemize}
		\item Objects are $\Omega_{Br}(n)=\{(p,A):p\in \Omega_{\tilde{+}}(n),\ A\in B_{len(p)}\}$ for $n\in \Z_{\geq0}$
		\item Symmetric group action $\sigma(p,A)=(\sigma p,A)$ for $(p,A)\in \Omega_{Br}(n)$ and $\sigma\in \Sigma_n$
		\item Composition given for $(p_n,A_n)\in \Omega_{Br}(n)$ and $(p_m,A_m)\in \Omega_{Br}(m)$ by
		\begin{equation}
			(p_n,A_n)\circ_i(p_m,A_m)=(p_n\circ_ip_m,A_n\circ_{t_{p_n}^{-1}(i)}A_m)
		\end{equation}
		\item Unit functors are given by $1\in \mathds{1}\mapsto (id_1,id_1)\in \Omega_{\tilde{+}}\times B_1$ viewed as an element of $\Omega_{Br}(1)$.  
	\end{itemize}
We view $p$ as a coloring of the strands of $A$ by the characters in $p$.  For $(p_n,\B_n)\in \Omega_{Br}(n)$ and $(p_m,\B_m)\in \Omega_{Br}(m)$, the composition $(p_n,A_n\circ_i(p_m,A_m)$ is then given by replacing the strand in $A_n$ labeled by $x_i$ with the braid $A_m$ made very thin, and coloring the resulting braid by $p_n\circ_ip_m\in \Omega_{\tilde{+}}(n+m-1)$.  
\end{definition}

\begin{proof*}
	It is clear by definition of length that $len(p_n\circ_i p_m)=len(p_n)+len(p_m)-1$ so this composition is well defined.  Then we need only check vertical and horizontal associativity, unit, and symmetric group equivariance.  Since this is structurally identical to the proof for definition \ref{wCoB}, the remainder of this proof is instead given below, after some notational remarks.
\end{proof*}

\begin{definition}
	In order to work with uniform notation for both cases, for $u_n\in \Omega_{\tilde{+}}(n)$, we will define 
	\begin{align}
		t_{u_n}=u_n&&len(u_n)=n
	\end{align}
\end{definition}

\begin{lemma}
	Suppose that $u_l\in \Omega_{\tilde{+}}(l)$ (resp. $\Omega_{br}(l)$), $u_m\in \Omega_{\tilde{+}}(m)$ (resp. $\Omega_{br}(m)$) and $u_n\in \Omega_{\tilde{+}}(n)$ (resp. $\Omega_{br}(n)$.  Then the following identities hold:
	\begin{itemize}
		\item Suppose $1\leq i\leq l$ and $1\leq j\leq m$.  Then 
		\begin{equation}\label{vertical1}
			t_{u_l\circ_iu_m}(t_{u_l}^{-1}(i)+t_{u_m}^{-1}(j)-1)=i+j-1
		\end{equation}
		since $t_{u_l}^{-1}(i)+t_{u_m}^{-1}(j)-1$ falls in the $i$-th block, at $t_{u_l}^{-1}(i)$.  
		\item Suppose $1\leq i<j\leq l$ and $t_{u_l}^{-1}(i)<t_{u_l}^{-1}(j)$.  Then 
		\begin{equation}\label{horizontal1a}
			t_{u_l\circ_ju_m}(t_{u_l}^{-1}(i))=i
		\end{equation}
		since $t_{u_l}^{-1}(i)$ falls left of the $j$-th block, which is between $t_{u_l}^{-1}(j)$ and $t_{u_l}^{-1}(j)+len(u_m)$.  Further, we have $t_{u_l}^{-1}(i)+len(u_m)<t_{u_l}^{-1}(j)+len(u_m)$ so 
		\begin{equation}\label{horizontal1b}
			t_{u_l\circ_iu_m}(t_{u_l}^{-1}(j)+len(u_m)-1)=j+len(u_m)-1
		\end{equation}
		since $t_{u_l}^{-1}(j)+len(u_m)-1$ falls right of the $j$-th block. 
		\item Suppose $1\leq i<j\leq l$ and $t_{u_l}^{-1}(i)>t_{u_l}^{-1}(j)$.  Then 
		\begin{equation}\label{horizontal2a}
			t_{u_l\circ_ju_n}(t_{u_l}^{-1}(i))=i+len(u_n)-1
		\end{equation}
		since $t_{u_l}^{-1}(i)$ falls right of the $j$-th block, which is between $t_{u_l}^{-1}(j)$ and $t_{u_l}^{-1}(j)+len(u_n)$.  Further, we have
		\begin{equation}\label{horizontal2b}
			t_{u_l\circ_i u_m}(t_{u_l}^{-1}(j))=j+len(u_m)-1
		\end{equation}
		since $t_{u_l}^{-1}(j)$ falls left of the $i$-th block, which is between $t_{u_l}^{-1}(i)$ and $t_{u_l}^{-1}(i)+len(u_m)$
	\end{itemize}
\end{lemma}

\begin{proof*}[Proof of \ref{wCoB} and \ref{omegabr}]  

Suppose that $(u_l,A_l)\in \wCoB_+(l)$ (resp. $\Omega_{br}(l)$), $(u_m,A_m)\in \wCoB_+(m)$ (resp. $\Omega_{br}(m)$) and $(u_n,A_n)\in \wCoB_+(n)$ (resp. $\Omega_{br}(n)$).

\textbf{vertical associativity}: Suppose that $1\leq i\leq l$ and $1\leq j\leq m$.  Then 
\begin{multline}
	((u_l,A_l)\circ_i(u_m,A_m))\circ_{i+j-1}(u_n,A_n)
	=(u_l\circ_iu_m,A_l\circ_{t_{u_l}^{-1}(i)}A_m)\circ_{i+j-1}(u_n,A_n)\\
	=((u_l\circ_i u_m)\circ_{i+j-1}u_n,(A_l\circ_{t_{u_l}^{-1}(i)}A_m)\circ_{t_{u_l\circ_iu_m}^{-1}(i+j-1)}A_n)\\
	=(u_l\circ_i(u_m\circ_ju_n),(A_l\circ_{t_{u_l}^{-1}(i)}A_m)\circ_{t_{u_l}^{-1}(i)+t_{u_m}^{-1}(j)-1}A_n)\\
	=(u_l\circ_i(u_m\circ_ju_n),A_l\circ_{t_{u_l}^{-1}(i)}(A_m\circ_{t_{u_m}^{-1}(j)}A_n))
	=(u_l,A_l)\circ_i(u_m\circ_ju_n,A_m\circ_{t_{u_m}^{-1}(j)}A_n)\\
	=(u_l,A_l)\circ_i((u_m,A_m)\circ_j(u_n,A_n))
\end{multline}
where, between the second and third lines, we use the identity presented in equation \ref{vertical1}. 

\textbf{horizontal associativity}: Suppose that $1\leq i<j\leq l$.  
\begin{multline}
	((u_l,A_l)\circ_i(u_m,A_m))\circ_{m+j-1}(u_n,A_n)
	=(u_l\circ_iu_m,A_l\circ_{t_{u_l}^{-1}(i)}A_m)\circ_{len(u_m)+j-1}(u_n,A_n)\\
	=((u_l\circ_iu_m)\circ_{len(u_m)+j-1}u_n,(A_l\circ_{t_{u_l}^{-1}(i)}A_m)\circ_{t_{u_l\circ_iu_m}^{-1}(len(u_m)+j-1)}A_n)
\end{multline}
\begin{multline}
	((u_l,A_l)\circ_j(u_n,A_n))\circ_i(u_m,A_m)
	=(u_l\circ_ju_n,A_l\circ_{t_{u_l}^{-1}(j)}A_n)\circ_i(u_m,A_m)\\
	=((u_l\circ_ju_n)\circ_iu_m,(A_l\circ_{t_{u_l}^{-1}(j)}A_n)\circ_{t_{u_l\circ_ju_n}^{-1}(i)}A_m)
\end{multline}
This equality for the $\Sigma$ (resp. $\Omega_{\tilde{+}}$) component follows from horizontal associativity in $\Sigma$ (resp. $\Omega_{\tilde{+}}$).  Then horizontal associativity reduces to proving equality of $(A_l\circ_{t_{u_l}^{-1}(i)}A_m)\circ_{t_{u_l\circ_iu_m}^{-1}(len(u_m)+j-1)}A_n$ and $(A_l\circ_{t_{u_l}^{-1}(j)}A_n)\circ_{t_{u_l\circ_ju_n}^{-1}(i)}A_m$.

\textbf{case: $t_{u_l}^{-1}(i)<t_{u_l}^{-1}(j)$}
\begin{multline}
	(A_l\circ_{t_{u_l}^{-1}(i)}A_m)\circ_{t_{u_l\circ_iu_m}^{-1}(len(u_m)+j-1)}A_n
	=(A_l\circ_{t_{u_l}^{-1}(i)}A_m)\circ_{t_{u_l}^{-1}(j)+len(u_m)-1}A_n\\
	=(A_l\circ_{t_{u_l}^{-1}(j)}A_n)\circ_{t_{u_l}^{-1}(i)}A_m
	=(A_l\circ_{t_{u_l}^{-1}(j)}A_n)\circ_{t_{u_l\circ_ju_n}^{-1}(i)}A_m
\end{multline}
where, in the first equality, we use the identity presented in equation \ref{horizontal1b}, between the first and second lines, we use horizontal associativity in $B$, and in the final equality, we use the identity presented in equation \ref{horizontal1a}. 

\textbf{case: $t_{u_l}^{-1}(i)>t_{u_l}^{-1}(j)$}
\begin{multline}
	(A_l\circ_{t_{u_l}^{-1}(i)}A_m)\circ_{t_{u_l\circ_iu_m}^{-1}(len(u_m)+j-1)}A_n
	=(A_l\circ_{t_{u_l}^{-1}(i)}A_m)\circ_{t_{u_l}^{-1}(j)}A_n\\
	=(A_l\circ_{t_{u_l}^{-1}(j)}A_n)\circ_{len(p_n)+t_{u_l}^{-1}(i)-1}A_m
	=(A_l\circ_{t_{u_l}^{-1}(j)}A_n)\circ_{t_{u_l\circ_ju_n}^{-1}(i)}A_m
\end{multline}
where, in the first equality, we use the identity presented in equation \ref{horizontal2b}, between the first and second lines, we use horizontal associativity in $B$, and in the final equality, we use the identity presented in equation \ref{horizontal2a}.

Then it follows that horizontal associativity holds in all cases. 

\textbf{symmetric group equivariance}: Suppose that $\sigma\in \Sigma_n$ and $\sigma_m\in \Sigma_m$.  Then 
\begin{multline}
	\sigma_n(u_n,A_n)\circ_i\sigma_m(u_m,A_m)=(\sigma_n u_n,A_n)\circ_i(\sigma_m u_m,A_m)\\
	=(\sigma_nu_n\circ_i\sigma_mu_m,A_n\circ_{t_{\sigma_n u_n}^{-1}(i)}A_m)=(\sigma_n\circ_i\sigma_n)(u_n\circ_{\sigma_n^{-1}(i)}u_m,A_n\circ_{t_{u_n}^{-1}\sigma_n^{-1}(i)}A_m)\\
	=(\sigma_n\circ_i\sigma_n)((u_n,A_n)\circ_i(u_m,A_m))
\end{multline}  

\textbf{unit}: It follows immediately from unit axiom on $\Sigma\times B$ (resp. $\Omega_{\tilde{+}}\times B$ that the unit axiom holds on objects in $\wCoB_+$ (resp. $\Omega_{br}$), since the unit and composition agrees with $\Sigma\times B$ (resp. $\Omega_{\tilde{+}}\times B$) on the level of objects. 

\end{proof*}

\begin{definition}
	For $(u,A)\in \Omega_{br}(n)$ and $1\leq i\leq len(u)$ define
	\begin{equation}
		\delta_i^p=\left\{\begin{array}{cc}
			*\in B_0&if\ t_u(i)=0\\
			id_1\in B_1& \text{otherwise}\end{array}\right.
	\end{equation}
	Define the \underline{underlying colored braid} of $(u,A)$ to be
	\begin{equation}
		\pi(u,A)=(\tau(u),\circ[A;\delta_1^u,...,\delta_{len(u)}^u])\in \wCoB_+(n)
	\end{equation}
	where $\circ[\_;\_]$ is the full composition in $B$.  One thinks of this map as removing the parentheses, and all strands colored by $0$.  
\end{definition}

\begin{proposition}
	The underlying braid map gives a morphism of symmetric operads in the category of groupoids 
	\begin{equation}
		\pi:\Omega_{Br}\to \wCoB_+
	\end{equation}
	where $\Omega_{Br}$ is considered as an operad in discrete groupoids.  
\end{proposition}

\begin{proof*}
It is clear that $\pi$ maps the unit in $\Omega_{Br}$ to the unit in $\wCoB_+$.  

Suppose that $(u_n,A_n)\in \Omega_{Br}(n)$ and $(u_m,A_n)\in \Omega_{Br}(m)$.  First note that since $\tau$ is a morphism of operads, we have $\tau(u_n)\circ_i\tau(u_m)=\tau(u_n\circ_iu_m)$.  We also make the following notes:
	\begin{itemize}
	\item If $1\leq j< t_{u_n}^{-1}(i)$, then $t_{u_n\circ_iu_m}(j)=t_{u_n}(j)$.  In particular, $\delta_{j}^{u_n\circ_iu_m}=\delta^{u_m}_j$.
	\item If $1\leq j\leq m$, then $t_{u_n\circ_i u_m}(t_{u_n}^{-1}(i)+j)=t_{u_m}(j)$.  In particular, $\delta^{u_n\circ_i u_m}_{t_{u_n}^{-1}(i)+j}=\delta^{u_m}_j$. 
	\item If $t_{u_n}^{-1}(i)< j\leq n$, then $t_{u_n\circ_iu_m}(len(u_m)+j-1)=t_{u_n}(j)$.  In particular, $\delta^{u_n\circ_i u_m}_{len(p_m)+j-1}a=\delta^{u_n}_j$.
\end{itemize}
Then, using this and associativity in $B$, we have 
\begin{multline}
	\circ[A_n;\delta^{u_n}_1,...,\delta^{u_n}_{len(n)}]\circ_{\tau(u_n)^{-1}(i)}\circ[A_m;\delta_1^{u_m},...,\delta^{u_m}_{len(m)}]\\
	=\circ[A_{t_{u_n}^{-1}(i)}A_m;\delta_1^{u_n},...,\delta^{u_n}_{t_{u_n}^{-1}(i)-1},\delta^{u_m}_1,...,\delta^{u_m}_{len(m)},\delta^{u_n}_{t_{u_n}^{-1}(i)+1},...,\delta^{u_n}_{len(n)}]\\
	=\circ[A_n\circ_{t_{u_n}^{-1}(i)}A_m;\delta^{u_n\circ_i u_m}_1,...,\delta^{u_n\circ_iu_m}_{len(u_n)+len(u_m)-1}]
\end{multline}
Then we conclude that $\pi(u_n,A_n)\circ_i\pi(u_m,A_m)=\pi((u_n,A_n)\circ_i(u_m,A_m))$, so $\pi$ is a morphism of operads on the level of objects.  Since $\Omega_{Br}$ is a discrete operad, there is nothing to check on the level of morphisms.  
\end{proof*}

\begin{definition}
	The \underline{weak braiding operad} is the operad in the category of groupoids defined by 
	\begin{equation}
		\wPaB_{\tilde{+}}=\pi^{-1}\wCoB_+
	\end{equation}
	We note that $\Omega_{Br}(n)$ and $\wPaB_{\tilde{+}}(n)$ share the same objects, and for $(p,A),(q,C)\in \wPaB_{\tilde{+}}(n)$, the morphism spaces are given by 
	\begin{equation}
		\Hom_{\wPaB_{\tilde{+}(n)}}(p,q)=Hom_{\wCoB_+(n)}(\pi(p,A),\pi(q,C))
	\end{equation}
\end{definition}

\begin{definition}
	A \underline{weak braided monoidal category} is an $\wPaB_{\tilde{+}}$ algebra.  
\end{definition} 

\begin{lemma}There is a morphism of symmetric operads
	\begin{equation}
		\wPaB_{\tilde{+}}\to\PaB_{\tilde{+}}
	\end{equation}
	induced by the commutative square 
	\begin{equation}
		\begin{tikzcd}
			\Omega_{Br}\arrow[r,"\pi"]\arrow[d,swap,"\Omega_f"]&\wCoB_{+}\arrow[d,"f"]\\
			\Omega_{\tilde{+}}\arrow[r,"\tau"]&\CoB_{+}
		\end{tikzcd}
	\end{equation}
	and functoriality of the limit, where $f$ and $\Omega_f$ are the morphisms of symmetric operads which forgets the braid component on objects.
\end{lemma}

\begin{corollary}
	If $\mathcal{C}$ is a braided monoidal category, then $\mathcal{C}$ is weak braided.
\end{corollary}
\begin{proof*}
	Since $\mathcal{C}$ is a braided monoidal category, there is morphism of symmetric operads $\phi:\PaB_{\tilde{+}}\to End_{\mathcal{C}}$.  Then we obtain morphism of symmetric operads 
	\begin{equation}\begin{tikzcd}
		\wPaB_{\tilde{+}}\arrow[r]&\PaB_{\tilde{+}}\arrow[r,"\phi"]&End_{\mathcal{C}}
	\end{tikzcd}\end{equation}
	so $\mathcal{C}$ obtains the structure of a $\wPaB_{\tilde{+}}$ algebra, and is a weak braided monoidal category. 
\end{proof*}

\begin{remark}
	Suppose that $\mathcal{C}$ is a monoidal category.  Then there is morphism of symmetric operads $\Omega_{\tilde{+}}\to End_{\mathcal{C}}$.  Suppose $p\in \Omega_{\tilde{+}}(n)$, and $X_1,....,X_n\in \mathcal{C}$.  Then we will denote 
	\begin{equation}
		p_\mathcal{C}(X_1,....,X_n)=\phi_n(p)(X_1,...,X_n)
	\end{equation}
	We think of this as tensoring the $X_i$ together according to the monomial defined by $p$, with the monoidal unit $\mathds{1}$ for each instance of $x_0$ in $p$.  Similarly, if $\mathcal{C}$ is weak-braided monoidal, so that there is morphism of symmetric operads $\phi:\wPaB_{\tilde{+}}\to End_{\mathcal{C}}$, and $(p,\B)\in \wPaB_{\tilde{+}}(n)$, then we will denote 
	\begin{equation}
		(p,\B)_\mathcal{C}(X_1,...,X_n)=\phi_n(p,\B)(X_1,...,X_n)
	\end{equation} 
\end{remark}

\section{Weak braiding for extensions}\label{weak braiding for extensions}

	\subsection{Module structure on coequilizers}
\begin{remark}\label{colimfunc}
Suppose that $I,\mathcal{C}$ are categories, and let $\mathcal{C}^I_0$ be a subcategory of the diagram category $\mathcal{C}^I$ on which the colimit of the functor $F\in \mathcal{C}^I_0$ exists.  Then the colimit is functorial.  That is, it induces a functor 
	\begin{equation}
		colim:\mathcal{C}^I_0\to \mathcal{C}
	\end{equation} 
	which takes a diagram to its colimit.  Notice that multiple such functors may exist, given by distinct (but uniquely isomorphic) choices of colimit for each diagram.

	Suppose that $I$ is a small category, and $\mathcal{C}^I_0$ is a subcategory of $\mathcal{C}^I$ for which colimits exist.  Suppose that $colim,colim':\mathcal{C}^I_0\to \mathcal{C}$ are functors as above. Then there is unique natural isomorphism $colim\to colim'$, induced by the universal property of the colimit. 
\end{remark}

\begin{definition} 
	For $n>0$, let $I_n$ be the two object category with $n$ morphisms $\mu^{(i)}:1\to 2$.  
	\begin{equation}\begin{tikzcd}[column sep=1.8cm]
			1\arrow[loop left]\arrow[r,shift left=0.3cm,"\mu^{(1)}"]\arrow[r,swap, shift left=-0.3cm,"\mu^{(n)}"]\arrow[r,shift left=-0.3cm]\arrow[r, draw=none,shift left=-0.25cm, "\vdots"]&2\arrow[loop right]
	\end{tikzcd}\end{equation}
Suppose that $\mathcal{C}$ is a category, and $F\in \mathcal{C}^{I_n}$.  Then a \underline{coequalizer} for the morphisms $F(\mu^{(i)}):F(1)\to F(2)$ is a colimit for the diagram $F$, which we denote by $coeq(F(\mu^{(i)}))_{i=1,...,n}$.
\end{definition}
	
	\begin{remark} 
		In the following we will assume that $\mathcal{C}$ is a braided monoidal category, and for any $X\in \mathcal{C}$, the functors $\cdot\boxtimes X$ and $X\boxtimes \cdot$ are right exact.
	\end{remark}
	\begin{definition}
		A \underline{commutative, associative algebra object} (or just algebra) in $\mathcal{C}$ consists of 
		\begin{itemize}
			\item $A\in \mathcal{C}$
			\item $\mu:A\boxtimes A\to A$
			\item $\iota_A:\mathds{1}\to A$
		\end{itemize}
		satisfying the following conditions
		\begin{enumerate}
			\item \textbf{Associativity}: \[\mu\circ (\mu\boxtimes 1_A)\circ \underline{\mathcal{A}}_{A,A,A}=\mu\circ (1_A\boxtimes \mu):A\boxtimes(A\boxtimes A)\to A\]
			\item \textbf{Commutativity}:
			\[\mu\circ \underline{\mathcal{R}}_A=\mu:A\boxtimes A\to A\]
			\item \textbf{Unit}:
			\[\mu\circ (\iota_A\boxtimes 1_A)\circ \underline{l}^{-1}_A=1_A:A\to A\]
		\end{enumerate}
	\end{definition}
	
	\begin{definition}
		Suppose $A\in \mathcal{C}$ is an algebra.  We define the category $\Rep_{\mathcal{C}}A$ or simply $\Rep A$ to be 
		\begin{itemize}
			\item \textbf{Objects}: $(X,\mu_X)$ with $X\in \mathcal{C}$, $\mu\in Hom_{\mathcal{C}}(A\boxtimes X,X)$ 
			satisyfing the following conditions
			\begin{enumerate}
				\item \textbf{Associativity}:
				\[\mu_X\circ (\mu\boxtimes 1_X)\circ \underline{\mathcal{A}}_{A,A,X}=\mu_X\circ (1_A\boxtimes \mu_X):A\boxtimes (A\boxtimes X)\to X\]
				\item \textbf{Unit}:
				\[\mu_X\circ (\iota_A\boxtimes 1_X)\circ \underline{l}_X^{-1}=1_X:X\to X\]
			\end{enumerate}
			\item \textbf{Morphisms}: For $(X,\mu_X),(Y,\mu_Y)\in \Rep A$ $Hom_{\Rep A}((X,\mu_X),(Y,\mu_Y))$ consists of $f\in Hom_\mathcal{C}(W_1,W_2)$ satisfying the following condition:
			\[f\circ \mu_{X}=\mu_{Y}\circ (1_X\boxtimes f)\]
		\end{itemize}
	\end{definition}

	\begin{proposition}\label{coeqRepA}
		Suppose that $(X,\mu^{(1)}_X),...,(X,\mu^{(n)}_X)\in \Rep A$, and $coeq(\mu^{(i)}_X)_{i=1,..,n}$ exists with coequilizer morphism $\eta$.  Then there is a unique $\mathcal{C}$ morphism 
		\begin{equation}
			\mu_{coeq}:A\boxtimes coeq(\mu^{(i)}_X)_{i=1,...,n}\to coeq(\mu^{(i)}_X)_{i=1,...,n}
		\end{equation}
		such that the following diagrams commute
		\begin{equation}\begin{tikzcd}
				A\boxtimes X\arrow[r,"\mu_X^{(i)}"]\arrow[d,swap,"1_A\boxtimes \eta"]&X\arrow[d,"\eta"]\\
				A\boxtimes coeq(\mu^{(i)}_X)_{i=1,...,n}\arrow[r,"\mu_{coeq}"]&coeq(\mu^{(i)}_X)_{i=1,...,n}
		\end{tikzcd}\end{equation}
		Moreover, $(coeq(\mu^{(i)}_X)_{i=1,...,n},\mu_{coeq})\in \Rep A$. 
		
	\end{proposition}
	
	\begin{proof*}
	Since $A\boxtimes\cdot$ is right exact, $(A\boxtimes (coeq(\mu^{(i)}_X)_{i=1,...,n},1_A\boxtimes \eta)$ is a coequilizer for the morphisms
		\begin{equation}
			1_A\boxtimes \mu^{(i)}_X:A\boxtimes(A\boxtimes X)\to A\boxtimes X
		\end{equation}
		We first claim that 
		\begin{equation}
			\eta\circ \mu^{(i)}_X\circ (1_A\boxtimes\mu^{(j)}_X)=\eta\circ \mu^{(i)}_X\circ (1_A\boxtimes\mu^{(k)}_X)
		\end{equation}
		For any $i,j\in \{1,...,n\}$.  Note that $\eta\circ \mu^{(i)}_X=\eta\circ \mu^{(j)}_X$, and, since $(X,\mu^{(j)}_X)$ is an $A$ module, we have $\mu^{(j)}_X\circ (1_A\boxtimes \mu_X)=\mu_X\circ (\mu\boxtimes 1_X)\circ \underline{\mathcal{A}}_{A,A,X}$. Then 
		\begin{multline}
			\eta\circ \mu^{(i)}_{X}\circ(1_X\boxtimes\mu^{(j)}_{X})=\eta\circ \mu^{(j)}_{X}\circ(1_X\boxtimes\mu^{(j)}_{X})=\eta\circ \mu^{(j)}_X\circ (\mu\boxtimes 1_X)\circ \underline{\mathcal{A}}_{A,A,X}\\
			=\eta\circ \mu^{(k)}_X\circ (\mu\boxtimes 1_X)\circ \underline{\mathcal{A}}_{A,A,X}
		\end{multline}
		Since $k\in \{1,...,n\}$ was arbitrary, we conclude that above claim.  Then by universal property of the coequilizer, there exists unique
		\begin{equation}
			\mu^{(i)}_{coeq}:A\boxtimes coeq(\mu^{(i)}_X)\to coeq(\mu^{(i)}_X)
		\end{equation}
		Such that $\mu^{(i)}_{coeq}\circ (1_A\boxtimes \eta)=\eta\circ \mu^{(i)}$.  But then $\mu^{(i)}_{coeq}\circ (1_A\boxtimes\eta)=\eta\circ \mu^{(i)}=\eta\circ \mu^{(j)}=\mu^{(j)}_{coeq}\circ (1_A\boxtimes \eta)$, so by uniqueness, we have $\mu^{(i)}_{coeq}=\mu^{(j)}_{coeq}$.  We shall define this common morphism to be $\mu_{coeq}$, and this makes all above diagrams commute.

		We need to confirm that $\mu_{coeq}$ satisfies the associativity and unit properties.  For the remainder of this proof, we will simply denote $coeq=coeq(\mu^{(i)}_X)_{i=1,...,n}$.  Then 
		
		\begin{multline}
			\mu_{coeq}\circ (\mu\boxtimes 1_{coeq})\circ \underline{\mathcal{A}}_{A,A,coeq}\circ (1_A\circ (1_A\circ \eta))=\mu_{coeq}\circ (\mu\boxtimes 1_{coeq})\circ ((1_A\boxtimes 1_A)\boxtimes \eta)\circ\underline{\mathcal{A}}_{A,A,X}\\
			=\mu_{coeq}\circ (1_A\boxtimes \eta)\circ (\mu\boxtimes 1_X)\circ \underline{\mathcal{A}}_{A,A,X}=\eta\circ \mu_X^{(i)}\circ  (\mu\boxtimes 1_X)\circ \underline{\mathcal{A}}_{A,A,X}\\
			=\eta \circ \mu^{(i)}_X\circ (1_A\boxtimes \mu^{(i)}_X)=\mu_{coeq}\circ (1_A\boxtimes \eta)\circ (1_A\boxtimes \mu^{(i)}_X)=\mu_{coeq}\circ (1_A\boxtimes \mu_{coeq})\circ (1_A\boxtimes(1_A\boxtimes \eta))
		\end{multline}
		By right exactness of $A\boxtimes\cdot$ and since $\eta$ is an epimorphism, $(1_A\boxtimes(1_A\boxtimes \eta))$ is also an epimorphism and we conclude that $\mu_{coeq}\circ (\mu\boxtimes 1_{coeq})\circ \underline{\mathcal{A}}_{A,A,coeq}=\mu_{coeq}\circ (1_A\boxtimes \mu_{coeq})$.  Then $\mu_{coeq}$ satisfies associativity. 
		
		Finally, to obtain the unit property: 
		\begin{multline}
			\mu_{coeq}\circ (\iota_A\boxtimes 1_{coeq})\circ l^{-1}_{coeq}\circ \eta=\mu_{coeq}\circ(\iota_A\boxtimes 1_X)\circ (1_{\mathds{1}}\boxtimes \eta)\circ l^{-1}_{X}\\
			=\mu_{coeq}\circ (1_A\boxtimes \eta)\circ (\iota_A\boxtimes 1_X)\circ l^{-1}_X=\eta\circ \mu^{(i)}_X\circ (\iota_A\boxtimes 1_X)\circ l^{-1}_X=\eta
		\end{multline}
		where we use the unit property of $\mu^{(i)}_X$.  Since $\eta$ is an epimorphism, we conclude that $\mu_{coeq}\circ (\iota_A\boxtimes 1_{coeq})\circ l^{-1}_{coeq}=1_{coeq}$, so $\mu_{coeq}$ satisfies the unit property.  Then $(coeq,\mu_{coeq})\in \Rep A$. 
		
	\end{proof*}
	
	\begin{remark}
		The previous definition of $\mu_{coeq}$ makes
		\begin{equation}
			\eta:X\to coeq(\mu^{(i)}_X)_{i=1,...,n}
		\end{equation}
		into a $\Rep A$ morphism. 
	\end{remark}

\subsection{Braid-tensor products}
Suppose that $\mathcal{C}$ is an abelian braided monoidal category.
\begin{definition}
		Suppose that $(X_1,\mu_{X_1}),...,(X_n,\mu_{X_{n}})\in \Rep A$ and $p\in \Omega_{\tilde{+}}(n)$.
		\begin{itemize}
			\item Define $X_0=A$.  
			\item Define 
			\begin{equation}
				p(A,X_1,...,X_n)=p^{len}_{\mathcal{C}}(\mathds{1},X_{t_p(1)},...,X_{t_p(len(p))})\end{equation}

			Note that this is $p(\mathds{1},X_{1},...,X_{n})$ with all instances of the monoidal unit $\mathds{1}$ replaced by $A$.  
			\item For $i=1,..,len(p)$, define the module map 
		\begin{equation}
			\mu^{(i)}_{p;X_1,...,X_n}:A\boxtimes p(A,X_1,...,X_n)\to p(A,X_1,...,X_n)
		\end{equation}
		By braiding over the first $i$-$1$ modules to use the $i$-th module map on $p(A,X_1,...,X_n)$
		\ctikzfig{munbraiddef}
		For each $i=1,...,len(p)$, $p(A,X_1,...,X_n)$ can be considered as an A module with module map $\mu^{(i)}_{p;X_1,...X_n}$.  
		\item Suppose $(p,\B)\in \wPaB_{\tilde{+}}(n)$.  Notice that $\B^{-1}$ acts as morphism
		\begin{multline}
			\B^{-1}:p(A,X_1,...,X_n)=p^{len}_{\mathcal{C}}(\mathds{1},X_{t_p(1)},...,X_{t_p(lenp)})\\
			\to p^{len}_{\mathcal{C}}(\mathds{1},X_{t_p\overline{\B}(1)},...,X_{t_p\overline{\B}(len(p))})
		\end{multline}
		For $i=1,...,len(p)$, we define the morphism
		\begin{multline}
			\mu^{(i)}_{p,\B;X_1,...,X_n}=\B\circ \mu^{(i)}_{p^{len};X_{t_p\overline{\B}(1)},...,X_{t_p\overline{\B}(len(p))}}\circ (1_A\boxtimes \B^{-1})\\
			A\boxtimes p(A,X_1,...,X_n)\to p(A,X_1,...,X_n)
		\end{multline}
		\end{itemize}
\end{definition}

\begin{remark}\label{mu12}
	Note that when $\mathcal{B}$ is the identity braiding on 2 elements, we have $\mu^{(i)}_{\nu,id;X_1,X_2}=\mu^{(i)}_{X_1,X_2}$ for $i=1,2$, agreeing with the notation of \cite{SVOAExtension} 
\end{remark}

\begin{lemma}\label{XYmod}
	Suppose that $(X_1,\mu_{X_1}),...,(X_n,\mu_{X_n})\in \Rep A$ and $(p,\mathcal{B})\in \wPaB_{\tilde{+}}(n)$.  Then $(p(A,X_1,...,X_n),\mu^{(i)}_{p,\mathcal{B};X_1,...,X_n})\in \Rep A$.
\end{lemma}

\begin{proof*} Note that by naturality of the braiding, this construction is independent of the parenthesization of $p$.  We begin by proving associativity. 
	
	\scalebox{0.9}{\parbox{\linewidth}{%
	\ctikzfig{munassoc}}}

	Second, the unit property. In the first step we make use of the naturality of the braiding and associators, and note that the $\mathds{1}$ strand may be brought through any other $X_j$ strand by using the triangle axiom and $\underline{l}_{X_j}\circ \underline{\mathcal{R}}_{X_j,A}\circ \underline{\mathcal{R}}_{A,X_j}=\underline{l}_{X_j}$.  In the second step we use the unit property of $\mu_{X_i}$.   
\scalebox{0.9}{\parbox{\linewidth}{%
	\ctikzfig{mununit}}}
\end{proof*}

\begin{definition}
	Suppose $l\geq0$.  Let $\Rep A^{I_l}_\mu$ be the full subcategory of $\Rep A^{I_l}$ such that for $F\in \Rep A^{I_l}_\mu$, 
	\begin{equation}
		F(\mu^{(i)})=\mu^{(j_i)}_{p,\B;X_1,...,X_n}
	\end{equation}
	For some $(p,\B)\in\PaB_{\tilde{+}}(n)$, $1\leq j_1,...,j_n\leq n$ and $X_1,...,X_n\in \Rep A$.  ie, $\Rep A^{I_l}_\mu$ consists of diagrams `built from' braid module maps. 
\end{definition}

\begin{lemma}\label{diagrammap}
	Suppose that $n\geq 0$ and $(p,\B)\in \wPaB_{\tilde{+}}(n)$.  The braid module maps define functors
	\begin{equation}
		\mu_{p,\B}:\Rep A^{\times n}\to \Rep A^{I_{len(p)}}_\mu
	\end{equation}
If $X_1,...,,X_n\in \Rep A$, then $\mu_{p,\B}$ is equivariant in the sense that for $\sigma\in S_n$, 
	\begin{equation}
		\mu_{\sigma(p,\B)}(X_1,...,X_n)=\mu_{(p,\B)}(X_{\sigma(1)},...,X_{\sigma(n)})
	\end{equation}
		
\end{lemma}

\begin{proof*}
	Suppose $(p,\B_n)\in \wPaB_{\tilde{+}}(n)$.  Define 
	\begin{equation}
		\mu_{p,\B}:\Rep A^{\times n}\to \Rep A^{I_{len(p)}}_\mu
	\end{equation}
	Such that for $X_1,...,X_n\in \Rep A$, $\mu_{p,\B}(X_1,...,X_n)$ is given by the following diagram:
	\begin{equation}\begin{tikzcd}A\boxtimes p(A,X_1,...,X_n)\arrow[r,shift left=0.3cm,"\mu^{(1)}_{p,\B}"]\arrow[r,swap, shift left=-0.3cm,"\mu_{p,\B}^{(len(p))}"]\arrow[r,shift left=-0.3cm]\arrow[r, draw=none,shift left=-0.25cm, "\vdots"]&p(A,X_1,...,X_n)
	\end{tikzcd}\end{equation}
	For $Y_1,...,Y_n\in \Rep A$ and $\Rep A$ morphisms $f_i:X_i\to Y_i$, we define a morphism $\mu_{p,\B}(f_1,...,f_n)$ of the diagrams by the following commutative diagram for $i=1,...,n$.  
	
	\begin{equation}\begin{tikzcd}
			A\boxtimes p(A,X_1,...,X_n)\arrow[d,swap,"1_A\boxtimes p(id_A{,}f_1{,}...{,}f_n)"]\arrow[r,"\mu^{(i)}_{p,\B;X}"]&p(A,X_1,...,X_n)\arrow[d,"p(id_A{,}f_1{,}...{,}f_n)"]\\
			A\boxtimes p(A,Y_1,...,Y_n)\arrow[r,"\mu^{(i)}_{p,\B;Y}"]&p(A,Y_1,...,Y_n)
	\end{tikzcd}\end{equation}
	This diagram commutes by naturality of the braiding and since each $\mu_{X_i}$ is a $\Rep A$ homomorphism.  Then this defines a natural transformation between these two diagrams.  Next we confirm that $\mu_{p,\B}$ defined in this way gives a functor.  
	
	It is clear from the definition that $\mu_{p,\B}(1_{X_1},...,1_{X_n})$ is the identity natural transformation on diagrams.  Suppose that $Z_1,...,Z_n\in \Rep A$ and $g_i:Y_i\to Z_i$ are $\Rep A$ morphisms.  Then we have the following commutative diagrams for $i=1,...,n$:
	\begin{equation}\begin{tikzcd}[row sep= 1.8cm]
			A\boxtimes p(A,X_1,...,X_n)\arrow[d,swap,"1_A\boxtimes p(id_A{,}f_1{,}...{,}f_n)"]\arrow[r,"\mu^{(i)}_{p,\B;X}"]&p(A,X_1,...,X_n)\arrow[d,"p(id_A{,}f_1{,}...{,}f_n)"]\\
			A\boxtimes p(A,Y_1,...,Y_n)\arrow[r,"\mu^{(i)}_{p,\B;Y}"]\arrow[d,swap,"1_A\boxtimes p^{A}(id_A{,}g_1{,}...{,}g_n)"]&p(A,Y_1,...,Y_n)\arrow[d,"p(id_A{,}g_1{,}...{,}g_n)"]\\
			A\boxtimes p(A,Z_1,...,Z_n)\arrow[r,"\mu^{(i)}_{p,\B;Z}"]&p(A,Z_1,...,Z_n)
	\end{tikzcd}\end{equation}
	In particular, since 
	\begin{equation}
		p(id_A,g_1,...,g_n)\circ p(id_A,f_1,...,f_n)=p(id_A,(g_1\circ f_1),..., (g_n\circ f_n))
	\end{equation}
	we have
	\begin{equation}
		\mu_{p,\B}((g_1,...,g_n)\circ(f_1,...,f_n))=\mu_{p,\B}(g_1,...,g_n)\circ \mu_{p,\B}(f_1,...,f_n)
	\end{equation}
	Then $\mu_{p,\B}$ does indeed define a functor. 
	
Symmetric group equivariance is clear from the definition of $\mu_{(p,\B)}$.
\end{proof*}

\begin{lemma}\label{diagramfunctornaturality}	
	Suppose that $(p,\B)\in \wPaB_{\tilde{+}}(n)$ and $g\in B_n$.  Define 
	\begin{equation}
	p'(x_0,x_1,...,x_n)=p^{len}(x_0,x_{t_p\overline{g}^{-1}(1)},...,x_{t_p\overline{g}^{-1}(len(p))})
	\end{equation}
	Then $g$ induces a natural transformation $\phi_n(p,\B)\to \phi_n(p',g\B)$.  
	
\end{lemma}
\begin{proof*}
	Suppose that $(X_1,\mu_{X_1}),...,(X_n,\mu_{X_n})\in \Rep A$.  Notice that $g$ acts as a morphism $p_\mathcal{C}(A,X_1,...,X_n)\to p'_{\mathcal{C}}(A,X_1,...,X_n)$.  By definition of the braid module maps, we have we have the following family of commutative diagrams, indexed by $i=1,...,len(p)$:
	\begin{equation}\label{braiddiagramnaturality}\begin{tikzcd}
			A\boxtimes p'_{\mathcal{C}}(A,X_1,...,X_n)\arrow[d,swap,"1_A\boxtimes g^{-1}"]\arrow[r,"\mu^{(i)}_{p',g\B;X}"]&p_{\mathcal{C}}(A,X_1,...,X_n)\arrow[d,"g^{-1}"]\\
			A\boxtimes p_{\mathcal{C}}(A,X_1,...,X_n)\arrow[r,"\mu^{(i)}_{p,\B;X}"]&p_{\mathcal{C}}(A,X_1,...,X_n)
	\end{tikzcd}\end{equation}
	Suppose that $(Y_1,\mu_{Y_1}),...,(Y_n,\mu_{Y_n})\in \Rep A$, and $f_j:(X_j,\mu_{X_j})\to (Y_j,\mu_{Y_j})$ is a $\Rep A$ morphism for $j=1,...,n$.  We have the following two families of commutative diagrams, which agree on their left and right edges, by naturality of $g^{-1}:p'_\mathcal{C}\to p_\mathcal{C}$.  In particular, we see that $g^{-1}$ induces a natural isomorphism $\phi_n(p',g\B)\to \phi_n(p,\B)$ as in the diagram \ref{braiddiagramnaturality}.
	
	\begin{equation}\begin{tikzcd}
			A\boxtimes p'_{\mathcal{C}}(A,X_1,...,X_n)\arrow[d,swap,"1_A\boxtimes g^{-1}"]\arrow[r,"\mu^{(i)}_{p',g\B;X}"]&p_{\mathcal{C}}(A,X_1,...,X_n)\arrow[d,"g^{-1}"]\\
			A\boxtimes p_{\mathcal{C}}(A,X_1,...,X_n)\arrow[r,"\mu^{(i)}_{p,\B;X}"]\ar[d,swap,"1_A\boxtimes {p(1_A,f_1,...,f_n)}"]&p_{\mathcal{C}}(A,X_1,...,X_n)\ar[d,"{p(1_A,X_1,...,X_n)}"]\\
			A\boxtimes p_{\mathcal{C}}(A,Y_1,...,Y_n)\ar[r,"\mu^{(i)}_{p,\B;Y}"]&p(A,Y_1,...,Y_n)
	\end{tikzcd}\end{equation}

\begin{equation}\begin{tikzcd}
		A\boxtimes p'_{\mathcal{C}}(A,X_1,...,X_n)\arrow[d,swap,"1_A\boxtimes {p'(1_A,f_1,...,f_n)}"]\arrow[r,"\mu^{(i)}_{p',g\B;X}"]&p_{\mathcal{C}}(A,X_1,...,X_n)\arrow[d,"{p'(1_A,f_1,...,f_n)}"]\\
		A\boxtimes p'_{\mathcal{C}}(A,Y_1,...,Y_n)\arrow[d,swap,"1_A\boxtimes g^{-1}"]\arrow[r,"\mu^{(i)}_{p',g\B;Y}"]&p_{\mathcal{C}}(A,Y_1,...,Y_n)\arrow[d,"g^{-1}"]\\
		A\boxtimes p_{\mathcal{C}}(A,Y_1,...,Y_n)\arrow[r,"\mu^{(i)}_{p,\B;Y}"]&p_{\mathcal{C}}(A,Y_1,...,Y_n)
\end{tikzcd}\end{equation}
\end{proof*}

\begin{lemma}
	Suppose that $r\geq0$.  The diagrams $F\in \Rep A^{I_r}_{\mu}$ have colimits.  For every planar binary tree $Y$ with $r$ ingoing edges, one assigns a coequalizer functor 
	\begin{equation}
		coeq_Y:\Rep A^{I_r}_\mu\to \Rep A
	\end{equation}
For any two binary trees $Y,Y'$ with $r$ ingoing edges, there is unique natural isomorphism $coeq_Y\to coeq_{Y'}$.
\end{lemma}
\begin{proof*}
	Recall that $\mathcal{C}$ is abelian, and so has cokernels.  In particular, for any two morphisms $f,g:M\to N$ in $\mathcal{C}$, one has binary coequalizers given by $coker(f-g)$.  
	
	Let $For:\Rep_{\mathcal{C}}A\to \mathcal{C}$ be the forgetful functor.  Then any $F\in \Rep A^{I_r}_\mu$ determines a diagram $For\circ F\in \mathcal{C}^{I_r}$.  Suppose that $Y$ is a rooted binary tree with $r$ ingoing edges.  One obtains a coequalizer functor $\widetilde{coeq}_{Y}:\mathcal{C}^{I_r}\to \mathcal{C}$ by taking binary coequalizers in the order determined by the binary tree $Y$.  
	
	By proposition \ref{coeqRepA}, there is unique morphism $\mu_{coeq_Y}:A\otimes \widetilde{coeq}_{Y}(F)\to coeq_Y(F)$ such that $(\widetilde{coeq}_Y(F),\mu_{coeq_Y(F)})\in \Rep_{\mathcal{C}}A$ is a coequalizer for $F$.  
	
	By functoriality of the colimit, this defines a functor $coeq_Y:\Rep A^{I_r}_\mu\to \Rep A$, and for any two binary trees $Y,Y'$ with $r$ ingoing edges, there is unique $\Rep A$-natural isomorphism $coeq_Y\to coeq_{Y'}$.  
\end{proof*}

\begin{theorem}\label{Extensionweakbraiding}
	Suppose that $\mathcal{C}$ is an abelian braided monoidal category and $A$ is an associative commutative algebra in $\mathcal{C}$. Then $\Rep_\mathcal{C}A$ is weak braided.  
\end{theorem}

\begin{proof*}
	We must construct a morphism of symmetric operads $\wPaB_{\tilde{+}}\to End_{\Rep_\mathcal{C}A}$.  Suppose that $n\geq 0$.  Suppose that $(p,\B)\in \wPaB_{\tilde{+}}$.  Let $Y_p$ by the binary tree with $len(p)$ ingoing edges associated to $p$, as in \cite{Fresse2017} 6.1.2. Define the functor 
	\begin{equation}
		\phi_n:\wPaB_{\tilde{+}}\to Hom_{Cat}(\Rep A^{\times n},\Rep A)
	\end{equation}
	on objects by
	\begin{equation}
		\phi_n(p,\B)=coeq_{Y_p}\circ \mu_{p,\B}:\Rep A^{\times n}\to \Rep A^{I_{len(p)}}_\mu\to \Rep A
	\end{equation}

	Suppose that there is $\wPaB_{\tilde{+}}(n)$ morphism $g:(p,\B_p)\to (q,\B_q)$.  We will now construct the desired natural isomorphism $\phi_n(p,\B_p)\to \phi_n(q,\B_q)$.  Let 
	\begin{equation}
		p'(x_0,x_1,...,x_n)=p_{\mathcal{C}}^{len}(x_0,x_{t_p\overline{\B_p}(1)},...,x_{t_p\overline{\B_p}(len(p))})
	\end{equation}

	so that $\B_p^{-1}$ and $\B_q^{-1}$ act as morphisms $\B_p^{-1}:p_\mathcal{C}(A,X_1,...,X_n)\to p'_\mathcal{C}(A,X_1,...,X_n)$.  As in lemma \ref{diagramfunctornaturality}, $\B_p^{-1}$ induces a natural isomorphism $\phi_n(p,\B_p)\to \phi_n(p',id_{len(p)})$ which `unwinds' the braiding.  Let $\lambda$ be the (unique) natural isomorphism  $\lambda:coeq_{Y_p}\to coeq_{Y_{p'}}$.  Then the horizontal composition of natural isomorphisms with  $\B^{-1}:\mu_{p,\B_p}\to \mu_{p',id_{len(p)}}$ defines a natural isomorphism
	\begin{equation}\phi_n(p,\B_p)=coeq_{Y_p}\circ \mu_{p,\B_p}\to coeq_{Y_{p'}}\circ \mu_{p',id_{len(p)}}=\phi_n(p',id_{len(p)})\end{equation}
	which intertwines the coequalizer morphisms.  Similarly, define 
	\begin{equation}
		q'(x_0,x_1,...,x_n)=q_{\mathcal{C}}^{len}(x_0,x_{t_q\overline{\B_q}(1)},...,x_{t_q\overline{\B_q}(len(p))})
	\end{equation}
	so that $\B_q^{-1}$ induces natural isomorphism $\phi_{n}(q,\B_q)\to \phi_{n}(q',id_{len(q)})$ which intertwines the respective coequalizer morphisms, as in the case of $(p,\B_p)$. 
	
	Note that $\circ[\B_p;\delta^p_1,...,\delta^p_{len(p)}]^{-1}$ acts as an isomorphism $(\tau_p,\circ[\B_p;\delta^p_1,...,\delta^p_{len(p)}])\to (\tau_{p'},id_n)$ in $\wCoB_+$, and so it acts as isomorphism $(p,\B_p)\to (p',id_{len(p)})$ in $\wPaB_{\tilde{+}}(n)$.  Similarly $\circ[\B_q;\delta^q_1,...,\delta^q_{len(q)}]^{-1}$ acts as an isomorphism $(q,\B_q)\to (q',id_{len(q)})$.

	This gives the following commutative diagram:
	\begin{equation}\begin{tikzcd}
			(p,\B_p)\arrow[r,"g"]\arrow[d,swap,"\circ{[\B_p;\delta^p_1,...,\delta^p_{len(p)}]}^{-1}"]&(q,\B_q)\arrow[d,"\circ{[\B_q;\delta^q_1,...,\delta^q_{len(q)}]}^{-1}"]\\
			(p',id_{len(p)})\arrow[r,dotted,"id"]&(q',id_{len(q)})
	\end{tikzcd}\end{equation}
Where we conclude that the bottom morphism must be labeled by identity braid since every composition of $id_1\in B_1$ with an identity braid is an identity braid. Then $\tau_p'=\tau_q'$.  Recalling that $\phi_n(p',id_{len(p)})$ and $\phi_n(q',id_{len(p)})$ correspond to the usual (un braid-twisted) monoidal structure on $\Rep_\mathcal{C}A$ by remark \ref{mu12},  there is unique $\mathcal{C}$-natural isomorphism $\phi_n(p',id_{len(p)})\to\phi_n(q',id_{len(p)})$ intertwining the coequilizer maps which is given by (compositions of) the associativity and unit morphisms, which we will denote here just by $\mathcal{A}$.  We denote by $\phi_n(g)$ the $\mathcal{C}$-natural isomorphism $\phi_n(p,\B_p)\to \phi_n(q,\B_q)$ making the following diagram commute
\begin{equation}\begin{tikzcd}
		\phi_n(p,\B_p)\arrow[r,dotted,"\phi_n(g)"]\arrow[d,"\B_p^{-1}"]&\phi_n(q,\B_q)\arrow[d,"\B_q^{-1}"]\\
		\phi_n(p',id_{len(p)})\arrow[r,"\mathcal{A}"]&\phi_n(q',id_{len(q)})
\end{tikzcd}\end{equation}
Since each of these natural isomorphisms intertwines coequilizer morphisms, so does their composition.  Using this $\mathcal{C}$-natural isomorphism, we see that $\phi_n(q,\B_q)$ gives coequilizers for the diagrams given by $\mu_{p,\B_p}$.  By proposition \ref{coeqRepA}, $\phi_n(g)$ is then a $\Rep_\mathcal{C}A$-natural isomorphism $\phi_n(p,\B_p)\to \phi_n(q,\B_q)$.   

It is immediate that $\phi_n$ defined in this way is a functor by uniqueness of the constructed morphisms.  We need only confirm that $\phi_n$ gives a morphism of symmetric operads. 
	
	Denote by $\mu_{(p_n,\B_n)\circ_i(p_m,\B_m)}(X_1,...,X_{n+m-1})|_{j=i,...,i+m}$ the diagram obtained from \\
	$\mu_{(p_n,\B_n)\circ_i(p_m,\B_m)}(X_1,...,X_{n+m-1})$ by considering only morphisms 
	\begin{equation*}\mu^{(i)}_{(p_n,\B_n)\circ_i(p_m,\B_m)},...,\mu^{(i+m)}_{(p_n,\B_n)\circ_i(p_m,\B_m)}\end{equation*}
	 Notice that 
	\begin{multline}coeq_{Y_{p_m}}(\mu_{(p_n,\B_n)\circ_i(p_m,\B_m)}(X_1,...,X_{n+m-1})|_{j=i,...,i+m})\\
		=\mu_{(p_n,\B_n)}(X_1,...,X_{i-1},\phi_m(X_{i+1},...,X_{i+m}),X_{i+m+1},...,X_{n+m-1}))
	\end{multline}		
	Then 
	\begin{multline}\phi_n(p_n,\B_n)\circ_i\phi_m(p_m,\B_m)(X_1,...,X_{n+m-1})\\
		=\phi_n(p_n,\B_n)(X_1,...,X_{i-1},\phi_m(p_m,\B_m)(X_{i},...,X_{i+m}),X_{i+m+1},...,X_{n+m-1})\\
		=coeq_{Y_{p_n}}(\mu_{p_n,\B_n}(X_1,...,X_{i-1},\phi_m(X_{i+1},...,X_{i+m}),X_{i+m+1},...,X_{n+m-1})))\\
		=\phi_{n+m-1}((p_n,\B_n)\circ_i(p_m,\B_m))(X_1,...,X_{n+m-1})
		\end{multline}
		Where in the final equality, we use the previous note and notice that $Y_{p_n\circ_ip_m}=Y_{p_n}\circ_iY_{p_m}$ so that $coeq_{p_n\circ_i p_m}$ is exactly the coequilizer functor constructed by first taking the coequalizer $coeq_{Y_{p_m}}|_{j=i,...,i+m}$ and then taking the coequilizer of remaining morphisms in the order dictated by $coeq_{Y_{p_n}}$.  Then the $\phi_n$ respect composition on objects, and in particular, composition intertwines the respective coequilizer morphisms, since the morphisms are constructed using universal property of the coequalizer.  
		
		Since we have already proven existence and uniqueness of the necessary morphisms, it is immediate that the $\phi_n$ respect composition of morphisms.  
		
Finally, we check symmetric group equivariance.  Suppose that $(p,\B)\in \wPaB_{\tilde{+}}(n)$ and $\sigma\in \Sigma_n$.  Notice that $Y_p=Y_{\sigma p}$.  Then
 \begin{multline}
 	\phi_n(\sigma(p,\B))(X_1,...,X_n)=coeq_{Y_{\sigma p}}\circ \mu_{\sigma(p,\B)}(X_1,...,X_n)\\
 	=coeq_{Y_p}\circ\mu_{(p,\B)}(X_{\sigma(1)},...,X_{\sigma(n)})=\phi_n(p,\B)(X_{\sigma(1)},...,X_{\sigma(n)})
 \end{multline}
Then the $\phi_n:\wPaB_{\tilde{+}}\to Hom_{Cat}(\Rep A^{\times n},\Rep A)=End_{\Rep A}(n)$ defines a morphism of symmetric operads, so $\Rep A$ is a weak-braided category.  
\end{proof*}
				
\section{Braided G-crossed categories}\label{GCrossed}
We can now turn our attention to the special case of braided $G$-crossed categories.  These should be weak braided, in the sense that they are `braided up to group action'.  See \cite{HFTandGCrossed} and \cite{EGNOTensor} for definitions. 

\begin{theorem}\label{Gcrossedcase}
	Let $G$ be a group and suppose that $\mathcal{D}$ is a braided $G$-crossed fusion category with additive structure functors.  Then $\mathcal{D}$ is weak-braided. 
\end{theorem}

\begin{proof*}
	We can write 
	\begin{equation}
		\mathcal{D}=\bigoplus_{g\in G}\mathcal{D}^g
	\end{equation}
	Let $F_g:\mathcal{D}\to \mathcal{D}$ be the `g-twisting functor' so that, for $X\in \mathcal{D}^g$, the braiding acts as 
	\begin{equation}
		\mathcal{R}_{XY}: X\otimes Y\to F_g(Y)\otimes X
	\end{equation}

	Suppose $(p,\B)\in \wPaB_{\tilde{+}}(n)$.  Suppose that $g_1,...,g_n\in G$ and $X_i\in \mathcal{D}^{g_i}$ for $i=1,..,n$.

	Notice that $\B$ acts as a morphism
	\begin{multline}
		\B:p^{len}_{\mathcal{D}}(x_o,x_{t_p\overline{\B}(1)},...,x_{t_p\overline{\B}(len(p))})(\mathds{1},X_1,...,X_n)\\
		\to p^{len}_{\mathcal{D}}(x_o,x_{t_p(1)},...,x_{t_p(len(p))})(\mathds{1},F_1(X_1),...,F_n(X_n))
	\end{multline}
	
	Where  \begin{equation}
		F_j=F_{g^{\epsilon_{j,k_j}}_{i_{j,k_j}}}...F_{g^{\epsilon_{j,1}}_{i_{j,1}}} 
		\end{equation} for some choice of $i_{j,l}\in \{1,..,n\}$ and $\epsilon_{j,k_j}=\pm1$.  It follows from the hexagon, pentagon and triangle axioms of braided G-crossed categories that the functors $F_j$ and this morphism depend only on the choice of braid $\B$ and twisted sectors $g_i$.  
	Define 
	\begin{multline}
		\phi_n(p,\B_p)(\mathds{1},X_1,...,X_n)=p_\mathcal{D}(\mathds{1},F_1(X_1),...,F_n(X_n))\\
		=p^{len}(x_o,x_{t_p(1)},...,x_{t_p(len(p))})_{\mathcal{D}}(\mathds{1},F_1(X_1),...,F_n(X_n))
	\end{multline}
	as a functor $\mathcal{D}^{g_1}\times...\times \mathcal{D}^{g_n}\to \mathcal{D}$ and extend as an additive functor $\mathcal{D}^{\times n}\to \mathcal{D}$.  Notice that, by this definition, $\B^{-1}$ induces a natural isomorphism $\phi_n(p,\B_p)\to \phi_n(p',id_{len(p)})$, since $\B$ acts as a natural transformation on $\mathcal{D}$.  
	
	Suppose that $\B\in Hom_{\wPaB_{\tilde{+}}(n)}((p,\B_p),(q,\B_q))$.  Write
	\begin{align}
		p'=p^{len}(x_o,x_{t_p\overline{\B_q}(1)},...,x_{t_p\overline{B_p}(len(p))})\\
		q'=q^{len}(x_o,x_{t_q\overline{\B_q}(1)},...,x_{t_q\overline{\B_q}(len(q))})
	\end{align}
Recall that $\circ[\B_p;\delta^p_1,...,\delta^p_{len(p)}]^{-1}$ acts as isomorphism $(p,\B_p)\to (p',id_{len(p)})$ in $\wPaB_{\tilde{+}}(n)$.  Similarly $\circ[\B_q;\delta^q_1,...,\delta^q_{len(q)}]^{-1}$ acts as an isomorphism $(q,\B_q)\to (q',id_{len(q)})$.  
	This gives the following commutative diagram:
	\begin{equation}\begin{tikzcd}
			(p,\B_p)\arrow[r,"\B"]\arrow[d,swap,"\circ{[\B_p;\delta^p_1,...,\delta^p_{len(p)}]}^{-1}"]&(q,\B_q)\arrow[d,"\circ{[\B_q;\delta^q_1,...,\delta^q_{len(q)}]}^{-1}"]\\
			(p',id_{len(p)})\arrow[r,dotted,"id"]&(q',id_{len(q)})
	\end{tikzcd}\end{equation}

Noticing that $\phi_n(p',id_{len(p)})$ and $\phi_n(q',id_{len(p)})$ correspond to the usual (untwisted) monoidal structure on $\mathcal{D}$,  there natural isomorphism $\phi_n(p',id_{len(p)})\to\phi_n(q',id_{len(p)})$, which is given by (compositions of) the associativity and unit morphisms, which we will denote here just by $\mathcal{A}$.  We define $\phi_n(\B)$ to be the natural isomorphism $\phi_n(p,\B_p)\to \phi_n(q,\B_q)$ making the following diagram commute
\begin{equation}\begin{tikzcd}
		\phi_n(p,\B_p)\arrow[r,dotted,"\phi_n(\B)"]\arrow[d,"\B_p^{-1}"]&\phi_n(q,\B_q)\arrow[d,"\B_q^{-1}"]\\
		\phi_n(p',id_{len(p)})\arrow[r,"\mathcal{A}"]&\phi_n(q',id_{len(q)})
\end{tikzcd}\end{equation}
It is clear that this defines a functor $\wPaB_{\tilde{+}}(n)\to End_{\mathcal{D}}(n)$, and also immediate from the definition that it intertwines the symmetric group action and respects units.  

Suppose that $(p_n,\B_n)\in \wPaB_{\tilde{+}}(n)$, $(p_m,\B_m)\in \wPaB_{\tilde{+}}(m)$.  We note that we have a commutative diagram of the following form in $End_{\mathcal{D}}(n+m-1)$, where we shall use the shorthand $p^{len}(x_{\sigma(j)})$ for $p^{len}(x_0,x_{\sigma(1)},...,x_{\sigma_{len(p)}})$. 
\begin{equation}\begin{tikzcd}[column sep= 11]
	p_n^{len}(x_{t_{p_n}\overline{\B_n}(j)})\circ_i p_m^{len}(x_{t_{p_m}\overline{\B_m}(j)})\arrow[r,equal]\ar[d,"\B_n\circ^{End}_{i}\B_m"]&(p_n\circ_ip_m)^{len}(x_{t_{p_n\circ_ip_m}\overline{\B_n\circ_{t_{p_n}^{-1}(i)}\B_m}(j)})\ar[d,"\B_n\circ_{t_{p_n}^{-1}(i)}\B_m"]\\
	\phi_n(p_n,\B_n)\circ_i\phi_m(p_m,\B_m)\ar[r,dotted]&\phi_{n+m-1}((p_n,\B_n)\circ_i(p_m,\B_m))
\end{tikzcd}\end{equation}
Where the operadic composition on the left is in the operad $End_{\mathcal{D}}$.  We then notice that the composition in the braid operad $\B_n\circ_{\tau_{p_n}^{-1}(i)}\B_m$ induces exactly the morphism $\B_n\circ_i^{End}\B_m$ up to associator or unit morphisms, so $\phi$ respects composition up to associator and unit morphisms, so the left and right morphisms are equal up to associators and unitors.  But it is clear from the definition that $\phi_n(p_n,\B_n)\circ_i\phi_m(p_m,\B_m)$ and $\phi_{n+m-1}((p_n,\B_n)\circ_i(p_m,\B_m))$ have the same arrangement of parentheses and units.  Then $\phi:\wPaB_{\tilde{+}}\to End_{\mathcal{D}}$ respects composition, so it is a morphism of operads and defines a weak braiding on $\mathcal{D}$.
\end{proof*}

\begin{remark}
	In the following, we let $\mathcal{R}$ be the generating braid on two strands.  
\end{remark}

\begin{theorem}\label{Gcrossedsufficient}
	Suppose that
	\begin{enumerate}
		\item $\mathcal{D}$ is linear, additive, semisimple and pivotal
		\item $\mathcal{D}$ is weak-braided with additive braid-tensor functors
		\item $S$ is a set of monoidal functors $F:\mathcal{D}\to \mathcal{D}$
		\begin{enumerate}
		\item For all $F\in S$ there exists simple $X\in \mathcal{D}$ such that 
			\begin{equation}
			(\nu,\mathcal{R})(\cdot,X)=(\nu,id)(F(\cdot),X)
		\end{equation}
		\item For all simple $X\in \mathcal{D}$, there is unique $F\in S$ such that 
		\begin{equation}
			(\nu,\mathcal{R})(\cdot,X)=(\nu,id)(F(\cdot),X)
		\end{equation}
		\item The tensor unit $\mathds{1}\in \mathcal{D}$ satisfies 
		\begin{equation}
			(\nu,\mathcal{R})(\cdot,\mathds{1})=(\nu,id)(id_{\mathcal{D}}(\cdot),\mathds{1})
		\end{equation}
		\end{enumerate}
	\end{enumerate}
	Then $S$ is a group and $\mathcal{D}$ is braided $S$-crossed.
\end{theorem}
\begin{proof*}
	It is immediate that $S$ has an associative multiplication given by composition of functors.  By condition 3c, we conclude that $id_{\mathcal{D}}\in S$ so $S$ forms a monoid.  
	
	For $F\in S$, define $\mathcal{D}^F$ to be the full subcategory of $\mathcal{D}$ on which 
	\begin{equation}
		(\nu,\mathcal{R})(\cdot,X)=(\nu,id)(F(\cdot),X)
	\end{equation}

	By condition 3a, each $\mathcal{D}^F$ is non-empty.

	Suppose that $X_i\in \mathcal{D}^F$  By condition 2, 
	\begin{equation}
		(\nu,\mathcal{R})(\cdot,\bigoplus_i X_i)=\bigoplus_i(\nu,\mathcal{R})(\cdot,X_i)
		=\bigoplus_i(\nu,id)(F(\cdot),X_i)=(\nu,id)(F(\cdot),\bigoplus_iX_i)
	\end{equation}
	Then the $\mathcal{D}^F$ are closed under direct sum, and $\mathcal{D}^F$ is contained inside the full subcategory of $\mathcal{D}$ generated by simples in $\mathcal{D}^F$.  
	
	Conversely, suppose that $X\in \mathcal{D}$.  Write $X=\bigoplus X_i$ for simple $X_i$.  Suppose that there exists $X_j\not \in \mathcal{D}^F$.  Using uniqueness in condition 3b, we see that 
	\begin{equation}
		(\nu,id)(F(\cdot),X_j)\neq (\nu,\mathcal{R})(\cdot,X_j)
	\end{equation} 
	These are exactly the $j$-th projection of the following two functors:
	\begin{equation}
		(\nu,\mathcal{R})(\cdot,X)=(\nu,\mathcal{R})(\cdot,\bigoplus_i X_i)=\bigoplus_i(\nu,\mathcal{R})(\cdot,X_i)
	\end{equation}
	\begin{equation}
		(\nu,id)(F(\cdot),X)	=\bigoplus_i(\nu,id)(F(\cdot),X_i)
	\end{equation}
Then \begin{equation}(\nu,\mathcal{R})(\cdot,X)\neq (\nu,id)(F(\cdot),X)\end{equation}
Then $X\not \in \mathcal{D}^F$.  We conclude that if $X$ is not a direct sum of simple $X_i\in \mathcal{D}^F$, then $X\not\in \mathcal{D}^F$.  
	
	Then $\mathcal{D}^F$ is the same as the full subcategory of $\mathcal{D}$ generated under direct sum by simple $X\in \mathcal{D}^F$.  We conclude that, for distinct $F$ and $G$ and $X\in \mathcal{D}^F$, $Y\in \mathcal{D}^G$, $Hom_\mathcal{D}(X,Y)=0$.  Then 
	\begin{equation}
		\mathcal{D}=\bigoplus_{F\in S}\mathcal{D}^F
	\end{equation}

	The following proofs will be guided by ribbon diagrams, although we must take care to deal with some subtlety as we are in the weak braided setting.  We will not justify a graphical calculus here,  and instead only give a translation to morphisms in $\mathcal{D}$, so that ribbon diagrams can be used to motivate a formal proof in terms of morphisms.  
	
	For strand $X\in \mathcal{D}^F$, we may replace 
	\begin{equation}
		(\nu,\mathcal{R})(\cdot,X)=(\nu,id)(F(\cdot),X)
	\end{equation}
	Then each braid morphism with $X\in \mathcal{D}^F$ acts as 
	\begin{equation}
		\mathcal{R}:(\nu,id)(X,\cdot)\to (\nu,id)(F(\cdot),X)
	\end{equation}
	Weak braiding satisfies naturality as usual.  Then, for $X\in \mathcal{D}^F$, we may use braiding as usual, with any strand $Z$ or coupon $f$ replaced by $F(Z)$ or $F(f)$, respectively, upon crossing under a strand labeled by $X$.  Associativity morphisms will be suppressed, using coherence of monoidal categories.
	
	Suppose that $W\in \mathcal{D}^F$ and $X\in \mathcal{D}^G$.  We claim that $W\otimes X=(\nu,id)(X,Y)\in \mathcal{D}^{FG}$.  Suppose that $Y,Z\in \mathcal{D}$ and $f:Y\to Z$.  This proof is guided by the following equality ribbon diagrams.
	
		\ctikzfig{Gcrossedproduct}

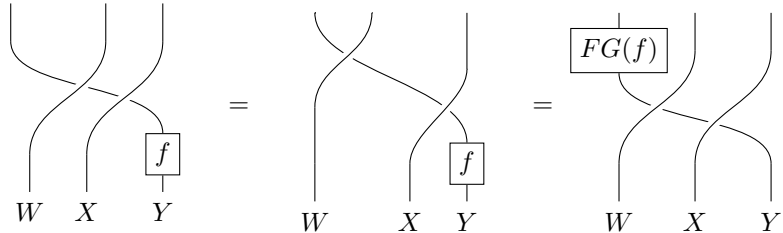
	\begin{multline}
		\mathcal{R}_{W\otimes X,Z}\circ (id_{W\otimes X}\otimes f)=(\mathcal{R}_{W,Z}\otimes id_X)(id_W\otimes \mathcal{R}_{X,Z}) (id_W\otimes id_X\otimes f)\\
		=(\mathcal{R}_{W,Z}\otimes id_X)(id_W\otimes G(f)\otimes id_Z)(id_W\otimes \mathcal{R}_{X,Z})\\
		=(FG(f)\otimes id_W\otimes id_X)(\mathcal{R}_{W,Z}\otimes id_X)(id_W\otimes \mathcal{R}_{X,Z})=(FG(f)\otimes id_{W\otimes X})\mathcal{R}_{W\otimes X,Z}
	\end{multline}
	Since $\mathcal{R}\circ_1id_2$ acts as morphism
	\begin{equation}(id,\nu(\nu(x_1,x_2),x_3))(X,Y,Z)\to (\mathcal{R}\circ_1id_2,\nu(x_3,\nu(x_1,x_2))(X,Y,Z)=(\mathcal{R},\nu)(X\otimes Y,Z)\end{equation}
	we conclude that $(\nu,\mathcal{R})(W\otimes X,\cdot)=(\nu,id)(FG(\cdot),W\otimes X)$. 
	
	We proceed by constructing inverses.  Suppose that $X\in \mathcal{D}^F$ is simple and consider its dual $X^*$.  $X^*$ is also simple, and so there is $G\in S$ such that 
	\begin{equation}
		(\nu,\mathcal{R})(\cdot,X^*)=(\nu,id)(G(\cdot),X^*)
	\end{equation}

	Suppose that $Y,Z\in \mathcal{D}$ and $f:Y\to Z$.  The construction is guided by the following ribbon diagram (omitting dimensional normalizations)
	
	\scalebox{0.9}{\parbox{\linewidth}{%
	\ctikzfig{Gcrossedinverse} }}
	
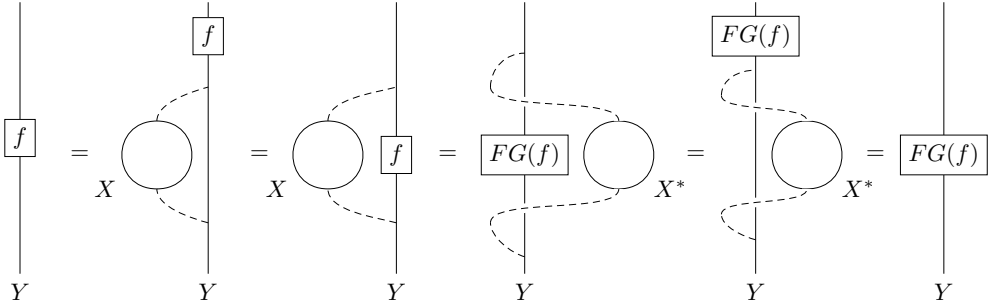
	\begin{multline}
		f=\frac{1}{dim(X)}fl_Y(ev_{X^*}\otimes id_Y)(coev_X\otimes id_Y)l^{-1}_Y\\
		=\frac{1}{dim(X)}l_Z(ev_{X^*}\otimes id_Y)(id_X\otimes id_{X^*}\otimes f)(coev_X\otimes id_Y)l^{-1}_Y\\
		=\frac{1}{dim(X)}l_Z(ev_{X^*}\otimes id_Y)\mathcal{R}_{X\otimes X^*,Y}^{-1}\mathcal{R}_{X\otimes X^*,Y}(id_X\otimes id_{X^*}\otimes f)(coev_X\otimes id_Y)l^{-1}_Y\\
		=\frac{1}{dim(X)}l_Z(ev_{X^*}\otimes id_Y)\mathcal{R}_{X\otimes X^*,Z}^{-1}(id_X\otimes id_{X^*}\otimes FG(f))\mathcal{R}_{X\otimes X^*,Y}(coev_X\otimes id_Y)l^{-1}_Y\\
		=\frac{1}{dim(X)}l_Z\mathcal{R}_{\mathds{1},Z}^{-1}(ev_{X^*}\otimes id_Y)(id_X\otimes id_{X^*}\otimes FG(f))(coev_X\otimes id_Y)\mathcal{R}_{\mathds{1},Z}l^{-1}_Y\\
		=l_Z\mathcal{R}_{\mathds{1},Z}^{-1}(id_\mathds{1}\otimes FG(f))\mathcal{R}_{\mathds{1},Y}l^{-1}_Y	=l_Z\mathcal{R}_{\mathds{1},Z}^{-1}\mathcal{R}_{\mathds{1},Y}(id_\mathds{1}\otimes FG(f))l^{-1}_Y\\
		=l_Z(id_\mathds{1}\otimes FG(f))l^{-1}_Y=FG(f)
	\end{multline}
	Since $Y,Z\in \mathcal{D}$ and $f:Y\to Z$ were arbitrary, we conclude that $GF=id_\mathcal{D}$.  We may similarly repeat this calculation, using instead $coev_{X^*}$ and $ev_X$, to obtain $GF=id_{\mathcal{D}}$.  Then each $F\in S$ has an inverse, so $S$ is a group.  Further, for simple $X\in \mathcal{D}^F$, we have $X^*\in \mathcal{D}^{F^{-1}}$.  
	
	Next, suppose that $X\in \mathcal{D}^F$ and $G\in S$. We claim that $G(X)\in \mathcal{D}^{GFG^{-1}}$.  Suppose that $Y,Z\in \mathcal{D}$, $f:Y\to Z$, and $W\in \mathcal{D}^G$.  This proof is guided by the following ribbon diagrams (with dimensional normalizations omitted).
	
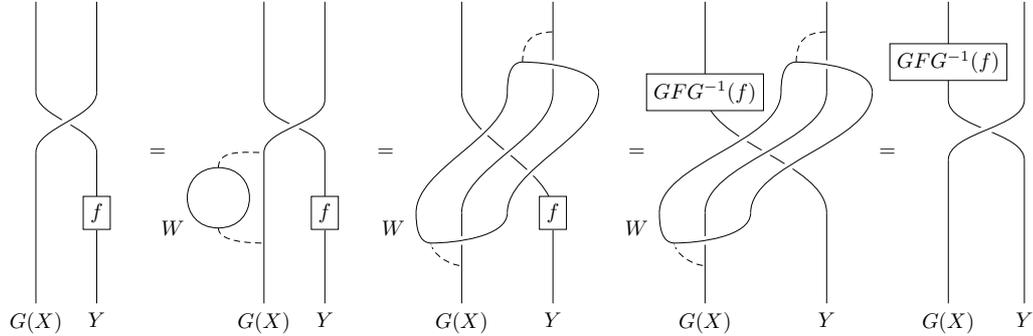
\scalebox{0.8}{\parbox{\linewidth}{%
	\ctikzfig{Gcrossedtarget}}}

We write this morphism instead as
	\begin{multline*}
		\mathcal{R}_{G(X),Y}(id_{G(X)}\otimes f)\\
		=\frac{1}{dim(W)}\mathcal{R}_{G(X),Z}(l_{G(X)}ev_{W^*}coev_{W}l^{-1}_{G(X)}\otimes id_Z)(1_{G(X)}\otimes f)\\
		=\frac{1}{dim(W)}\mathcal{R}_{G(X),Z}(l_{G(X)}\otimes id_Z)(ev_{W^*}\otimes id_{G(X)}\otimes id_{Z})\circ\\
		(id_{W}\otimes \mathcal{R}^{-1}_{W^*,G(X)}\mathcal{R}_{W^*,G(X)}\otimes id_{Z})(coev_{W}\otimes id_Z)(l^{-1}_{G(X)}\otimes id_Z)
		(id_{G(X)}\otimes f)
	\end{multline*}
	\begin{multline*}
		=\frac{1}{dim(W)}(id\otimes l_{G(X)})(id\otimes ev_{W^*}\otimes id_{G(X)})(id\otimes id_{W}\otimes \mathcal{R}^{-1}_{W^*,G(X)})\circ\\
		\mathcal{R}_{W\otimes G(X)\otimes W^*,Z}(id_{W}\otimes \mathcal{R}_{W^*,G(X)}\otimes id_{Z})(coev_{W}\otimes id_Z)(l^{-1}_{G(X)}\otimes id_Z)(id_{G(X)}\otimes f)
	\end{multline*}
	\begin{multline}
		=\frac{1}{dim(W)}(GFG^{-1}(f)\otimes id_{G(X)})(id_{GFG^{-1}(Y)}\otimes l_{G(X)})(id_{GFG^{-1}(Y)}\otimes ev_{W^*}\otimes id_{G(X)})\circ\\
		(id_{GFG^{-1}(Y)}\otimes id_{W}\otimes \mathcal{R}^{-1}_{W^*,G(X)})\mathcal{R}_{W\otimes G(X)\otimes W^*,Y}(id_{W}\otimes \mathcal{R}_{W^*,G(X)}\otimes id_{Z})\circ\\
		(coev_{W}\otimes id_Y)(l^{-1}_{G(X)}\otimes id_Y)
		=(GFG^{-1}(f)\otimes id_{G(X)})\mathcal{R}_{G(X),Y}
	\end{multline}
	Then we see that $(\nu,\mathcal{R})(G(X),\cdot)=(\nu,id)(GFG^{-1}(\cdot), G(X))$ so $G(X)\in \mathcal{D}^{GFG^{-1}}$.

	To see that the hexagon axiom is satisfied, notice that $\mathcal{R}\circ_1 id_2=(\mathcal{R}\otimes id)(id\otimes \mathcal{R})$ are the same braid, and so correspond to the same morphism (between appropriately chosen source and target space) under the $\wPaB_{\tilde{+}}$ algebra structure.  
\end{proof*}

\printbibliography 

\end{document}